\documentclass[final]{siamltex}

\usepackage{amsmath,amssymb}
\usepackage{booktabs}
\usepackage{graphicx}
\usepackage{listings}
\usepackage{algorithm}
\usepackage{algpseudocode} 
\usepackage{comment}
\usepackage{color}

\definecolor{dkgreen}{rgb}{0,0.6,0}
\definecolor{gray}{rgb}{0.5,0.5,0.5}
\definecolor{mauve}{rgb}{0.58,0,0.82}
\definecolor{pink}{rgb}{0.9,0,0.5}

\usepackage[
colorlinks=true,
citecolor=blue,
urlcolor=blue,
linktocpage,
pagebackref
]{hyperref}
\usepackage[all]{hypcap}

\newcommand{\code}[1]{\texttt{\small\color{black} #1}}

\newcommand{\Hcurl}{\ensuremath{H(\mathrm{curl})}}
\newcommand{\Hdiv}{\ensuremath{H(\mathrm{div})}}

\title{Non-Conforming Mesh Refinement for High-Order Finite Elements%
  \thanks{Performed under the auspices of the U.S. Department of Energy under
    Contract DE-AC52-07NA27344 (LLNL-JRNL-751849). The first author has been
    supported by the Ministry of Education, Youth and Sports of the Czech
    Republic under the RICE New Technologies and Concepts for Smart Industrial
    Systems, project No. LO1607.}}

\author{
Jakub \v{C}erven\'{y}
  \thanks{Center for Applied Scientific Computing (CASC),
    Lawrence Livermore National Laboratory, Livermore, CA 94551,
    \texttt{cerveny1@llnl.gov}, \texttt{dobrev1@llnl.gov},
    \texttt{tzanio@llnl.gov}}
\and Veselin Dobrev\footnotemark[2] %
\and Tzanio Kolev\footnotemark[2] %
}

\begin{document}
\maketitle

\begin{abstract}
  We propose a general algorithm for non-conforming adaptive mesh refinement
  (AMR) of unstructured meshes in high-order finite element codes. Our focus is
  on $h$-refinement with a fixed polynomial order.
  The algorithm handles triangular, quadrilateral, hexahedral and prismatic meshes
  of arbitrarily high order curvature, for any order finite element space in the
  de Rham sequence. We present a flexible data structure for meshes with hanging
  nodes and a general procedure to construct the conforming interpolation operator,
  both in serial and in parallel.
  The algorithm and data structure allow anisotropic refinement of tensor product elements
  in 2D and 3D, and support unlimited refinement ratios of adjacent elements.
  We report numerical experiments verifying the correctness
  of the algorithms, and perform a parallel scaling study to show that we can
  adapt meshes containing billions of elements and run efficiently on 393,000
  parallel tasks. Finally, we illustrate the integration of dynamic AMR into
  a~high-order Lagrangian hydrodynamics solver.
\end{abstract}



\pagestyle{myheadings}
\thispagestyle{plain}
\markboth{}{}

\section{Introduction} \label{sec:intro}

High-order finite element methods (FEM) are becoming increasingly important in the
field of scientific computing due to their increased accuracy and performance efficiency
on modern computer architectures \cite{ceed}. Many high-order applications are
interested in parallel adaptive mesh refinement (AMR) on unstructured meshes
containing quadrilateral and hexahedral elements. Such
elements are attractive for their tensor product structure
(enabling efficiency) and for their refinement flexibility (enabling {\em
anisotropic} refinement). However, as opposed to bisection-based methods for
local refinement of simplices, {\em hanging} nodes that occur after local
refinement of quadrilaterals and hexahedra are not easily avoided by further
refinement \cite{Schneiders98}. We are thus interested in handling
{\em non-conforming} (irregular) meshes \cite{Demkowicz89, Schonfeld95,
Heuveline07}, in which adjacent elements need not share a complete face or edge
and where some finite element degrees of freedom (DOFs) have to be constrained
to obtain a conforming, high-order $H^1$, $\Hcurl$ or $\Hdiv$ solution.

In this paper we present a set of software abstractions and algorithms for
handling large parallel non-conforming meshes and high-order function spaces,
that we believe strikes the right
balance between efficiency, ease of implementation and applicability to a wide
range of problems. In particular, with a single general algorithm we handle a
large class of 2D and 3D AMR meshes that consist of quadrilaterals and hexahedra,
possibly mixed with triangles and triangular prisms. Generally, the
elements can be refined by bisecting their edges, and we allow
anisotropic refinement of tensor product elements, meaning that e.g. a~
hexahedron can be split into 2 or 4 hexahedra only. In addition, the method we
propose is on a general discretization level, independent of the physics
simulation, and supports the whole de Rham sequence of finite element spaces, at
arbitrarily high order. Our approach is also compatible with curvilinear meshes,
as well as with high-order finite element techniques such as static condensation
and hybridization \cite{Dobrev17}. While the algorithms we present can be extended to
general $hp$-refinement, in this paper we focus exclusively on non-conforming
$h$-refinement with a fixed order $p \geq 1$.

The methodology of constraining hanging nodes in non-conforming meshes is well
established in the literature. Most authors opt to construct algebraic
constraints that express constrained degrees of freedom as linear combinations
of {\em true} degrees of freedom. The constraints are then applied to either the local
element matrices as part of the global matrix assembly \cite{Ainsworth97,
Demkowicz89, Carey76, libMesh06, Bank83}, after global assembly to ``condense''
the constrained DOFs out of the linear system \cite{dealII07}, or, more rarely,
at the shape function level \cite{Solin08} or the linear solver level
\cite{Suttmeier07}. Yet another possibility is to integrate the interpolation of
hanging nodes within a multigrid framework \cite{Brandt1977, Reps2017}. Our
approach is closest to the method of adjusting the global linear system, but
we cast the elimination of constrained DOFs as a form of variational restriction.
The practical advantage is that the stiffness matrix assembly stays unchanged and
completely independent of the treatment of AMR, and we also avoid intermixing AMR
with the numerical method of the simulation.

A large part of this paper is concerned with an efficient algorithm to construct
a~global operator that interpolates the constrained degrees of freedom from the
true degrees of freedom. The algorithm is based on composing local interpolation
matrices into a global one. It is easy to implement, works for a variety of finite
element spaces and element geometries and is unique in being able to handle 3D anisotropic
refinements as well as arbitrarily irregular meshes where adjacent elements are not
limited to a 2:1 refinement ratio. The resulting interpolation matrix is also
conveniently used to constrain duplicate nodes in parallel domain decomposition.

A modern
finite element solver needs to scale to hundreds of thousands of parallel
tasks. To achieve that, the mesh has to be fully distributed, as the global
mesh no longer fits on each processor. Historically, this is not been the
case in popular codes such as libMesh \cite{libMesh06}. The deal.II project
\cite{dealII07} has switched to the p4est library to distribute the mesh in 2011.
In addition to mesh distribution, an AMR solver at this level of parallelism
requires an extremely
efficient partitioning and load balancing algorithm that can be applied often.
High quality spectral type partitioners, such as ParMETIS \cite{Karypis98}, may
be prohibitively expensive for this purpose. We thus take inspiration from
octree-based approaches \cite{Flaherty97, p4est11, IBW2015} where only the coarse
elements (the octree roots) are shared by all parallel tasks and where the
refinement tree provides a natural way to enumerate elements. The resulting
sequence is equipartitioned and distributed among processors, which has been
shown to be a scalable way to load balance the mesh \cite{p4est11}. In contrast
to octree methods, however, we allow for more general refinement and support
additional types of elements. Our partitioning algorithm is thus closest to the
REFTREE method \cite{Mitchell07}.

The rest of the paper is organized as follows: in Section \ref{sec:ncmeshes}, we
define the class of meshes we support. In Section \ref{sec:amr} we describe the
general procedure to construct the AMR system through variational restriction. Section
\ref{sec:serialP} is devoted to the algorithm for constructing the AMR
interpolation operator, in serial. The internal data structures for
element-based AMR in the context of high-order finite elements are presented in
Section \ref{sec:data-serial}. Extending these methods to parallel settings is
the subject of Section \ref{sec:parallel}. Finally, Section \ref{sec:numerical}
presents numerical experiments that illustrate the performance of our
algorithms in practice.

\section{Non-conforming meshes} \label{sec:ncmeshes}

For the purposes of this paper, and following the definitions in
\cite{Mitchell05}, an \emph{element} $K$ is a (closed) triangle, quadrilateral,
triangular prism or hexahedron. An element contains $vertices$, $edges$, with
the usual definitions, and three dimensional elements contain $faces$. A mesh
$\{K_i\}_{i=1}^M$ is the union of $M$ elements, not necessarily all of the same
type, such that $\Omega = \cup_i K_i$ is a connected, bounded region in
$\mathbb{R}^2$ or $\mathbb{R}^3$.  Let $\langle K \rangle$ and $\partial K$
denote the interior and boundary of element $K$. We assume that the element
interiors do not intersect, i.e. $\langle K_i \rangle \cap \langle K_j \rangle =
\emptyset$, for $i \neq j$. A mesh is said to be \emph{conforming} if the set
$\partial K_i \cap \partial K_j$, $i \neq j$, is a common vertex, common edge,
common face or an empty set. A vertex of an element is called a \emph{hanging node}, if
it lies in the interior of an edge or face of another element. A mesh that
contains at least one hanging node is called \emph{non-conforming} or irregular.

To \emph{refine} an element $K_i$ means to replace it with at least two smaller
elements $K_{ik}$ of the same type, such that $K_i = \cup_k K_{ik}$. We call the
elements $K_{ik}$ \emph{child elements} and $K_i$ the \emph{parent element}. The
parent element is removed from the mesh after refinement, but it is stored
along with links to the child elements in a \emph{refinement tree}. In an
inverse, \emph{coarsening} process\footnote{A better term would in fact be \emph{derefinement},
meaning that only previously introduced refinements can be removed.}, the parent
element can later be reintroduced to the mesh, once its child elements are removed.

\begin{figure}
\begin{center}
  \includegraphics[width=0.31\textwidth]{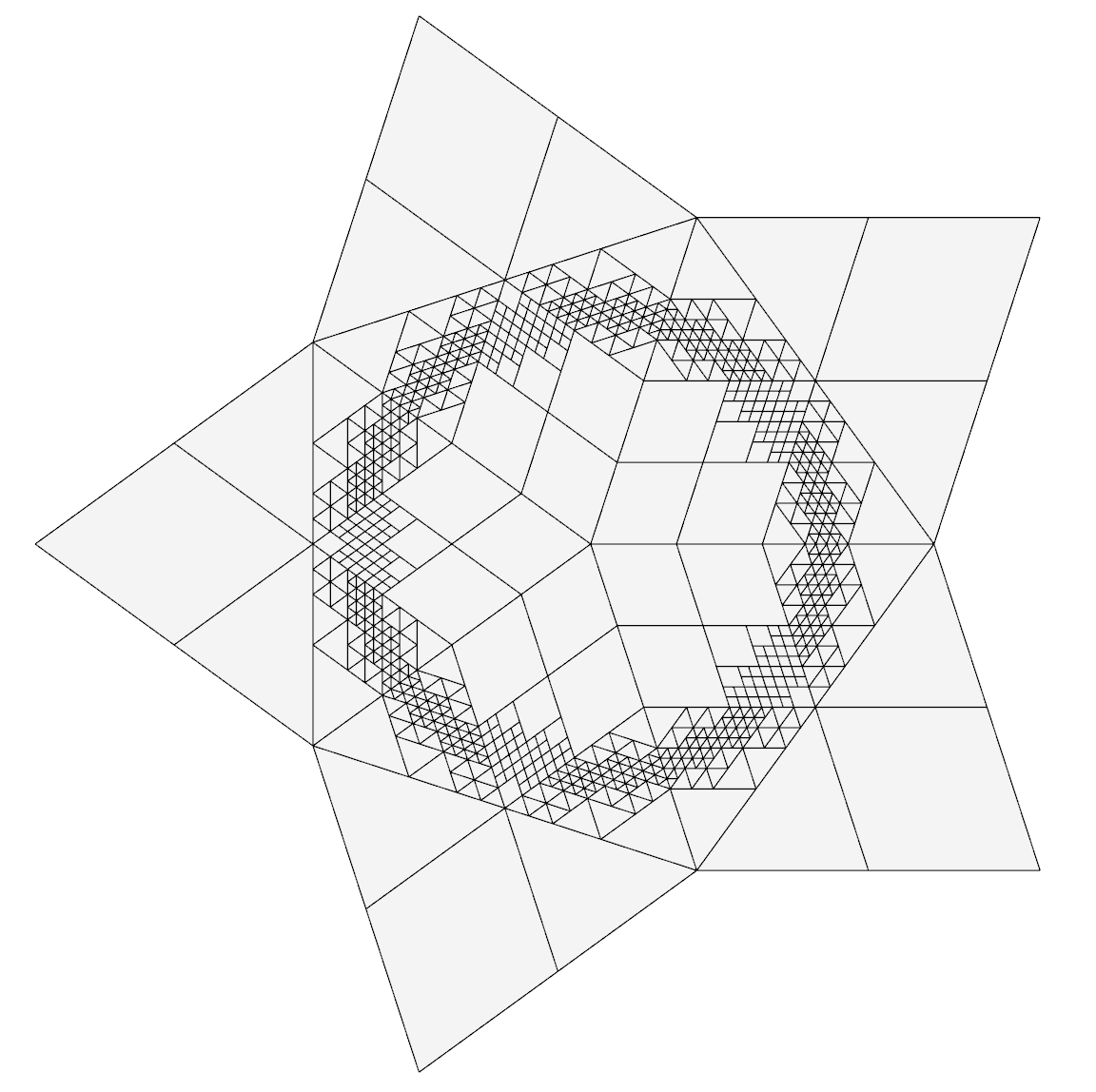} \hfill
  \includegraphics[width=0.31\textwidth]{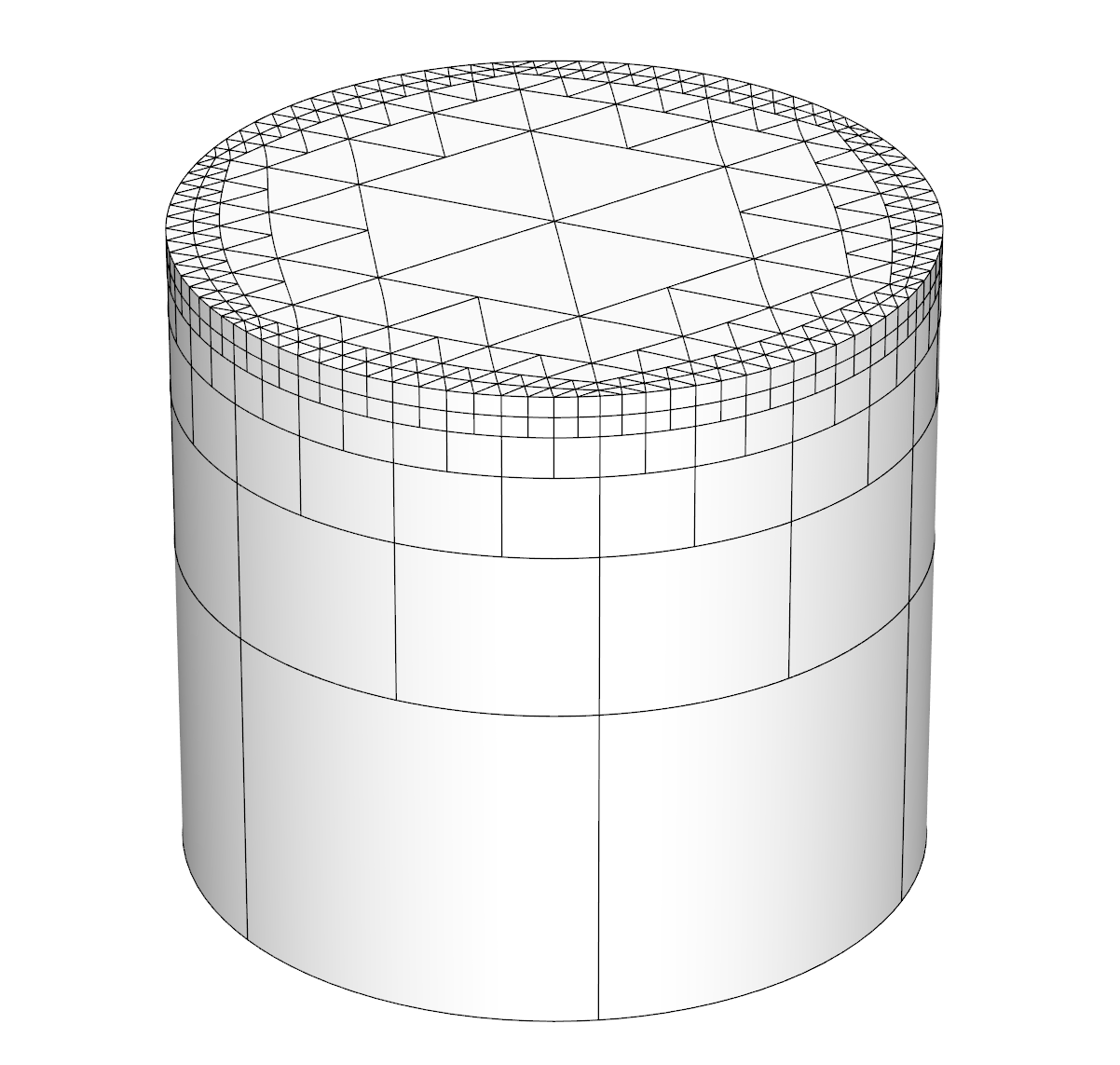} \hfill
  \includegraphics[width=0.31\textwidth]{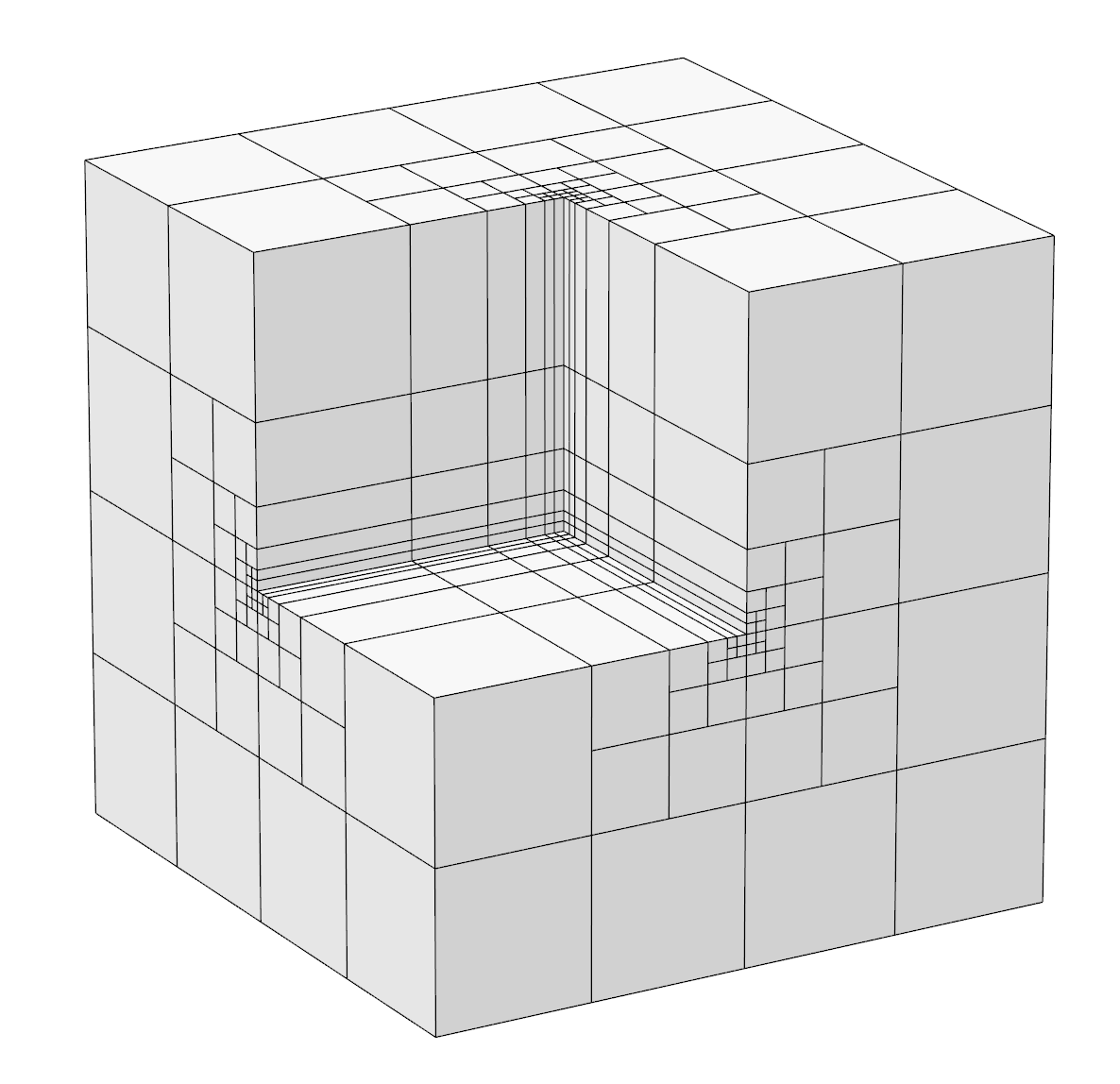}
\end{center}
\caption{Examples of non-conforming meshes handled by our algorithm:
  a mixed mesh containing triangles and quadrilaterals (left), a prismatic mesh
  obtained by refining 12 curved triangular prisms (middle), a hexahedral mesh
  after anisotropic refinement (right).}
\label{fig:ncmesh}
\end{figure}

When refining an element, we require that new vertices be placed in the center
of edges or faces of the parent element. We thus only support bisection of edges
during refinement. In this way a triangle can be split into 4 children, a
quadrilateral into 2 or 4 children, and a prism or hexahedron into 2, 4 or 8
children. We refer to refinements that split all edges of the parent element as
\emph{isotropic} refinement, and the remaining refinement options as
\emph{anisotropic}. Some examples of the types of non-conforming meshes we support
(obtained by refining initially conforming meshes) are shown in Figure \ref{fig:ncmesh}.

A non-conforming mesh is said to be \emph{consistent} if lower dimensional mesh
entities (i.e., vertices, edges, faces) that intersect are either identical, or
one is a proper subset of the other (see Figure \ref{fig:aniso} for an example
of an inconsistent mesh). We refer to the smaller entity as a \emph{slave} and
the larger containing entity as a~\emph{master}.  Lower dimensional entities of a
consistent mesh that are neither masters nor slaves are called \emph{conforming}
entities.

Note that a slave face has exactly one master face, but a 3D edge can have
multiple master edges. In that case we always take the smallest master edge,
which in turn may be a slave to its smallest master. This hierarchy will end
with a maximal master edge that is a not a slave itself. Note that any master
edge can be shared by several elements in a ``conforming'' manner. We call the
interface between slave faces and their master, or slave edges and their maximal
master the {\em coarse-fine interface}.

\section{Non-conforming AMR by variational restriction} \label{sec:amr}

In the continuous Ga\-ler\-kin method \cite{Ciarlet78}, one starts with a weak
variational problem: Find $u \in V$, such that
\begin{equation}
  a(u,v) = l(v) \qquad \forall v \in V \,, \label{eq:weak}
\end{equation}
where $a(\cdot,\cdot)$ and $f(\cdot)$ are bilinear and linear forms,
respectively, on the inner-product space $V$, typically a Sobolev space on the
domain $\Omega$ with imposed essential boundary conditions. The problem \eqref{eq:weak}
is discretized by constructing an approximate finite dimensional subspace $V_h
\subset V$, $\mathrm{dim}(V_h) = N$, on the mesh $\{K_i\}_{i=1}^M$ and
by assembling the stiffness matrix $A$ and load vector $f$, so that
\begin{equation}
  v_h^T A u_h = v_h^T f \qquad \forall v_h \in \mathbb{R}^{N}. \label{discreteprob}
\end{equation}
This is equivalent to the linear system $A u_h = f$, which yields the approximate solution vector $u_h$.

When constructing the space $V_h$, care must be taken to ensure that conformity
requirements of the space $V$ are met, so that $V_h$ is a subspace of~$V$. For
example, if $V$ is the Sobolev space $H^1$, the functions in $V_h$ must be kept
continuous across element boundaries. If $V$ is an $\Hcurl$ space, the
tangential component of the finite element vector fields in $V_h$ needs to be
continuous across element faces, etc. In high-order FEM this is achieved by
matching suitable \emph{shape functions} at degrees of freedom associated with
the mesh vertices, edges and faces, which together with DOFs in the element
interiors form the $N$ basis functions of $V_h$ \cite{SolinChapter4}.

When the mesh is non-conforming, the shape functions cannot be matched directly
at the coarse-fine interfaces. If we assign
degrees of freedom in the usual way on the conforming entities (including
vertices, and maximal master edges in 3D) and disregard any conformity in the interior
coarse-fine interfaces, we obtain a larger space $\widehat{V}_h \supset V_h$
that is ``broken'' along these interfaces, i.e., slave entities have DOFs that are
independent of their master DOFs. We call $\widehat{V}_h$ the \emph{partially
conforming space}, with $\dim(\widehat{V}_h) = \widehat{N} > N$.

\begin{figure}
\begin{center}
  \includegraphics[width=0.3\textwidth]{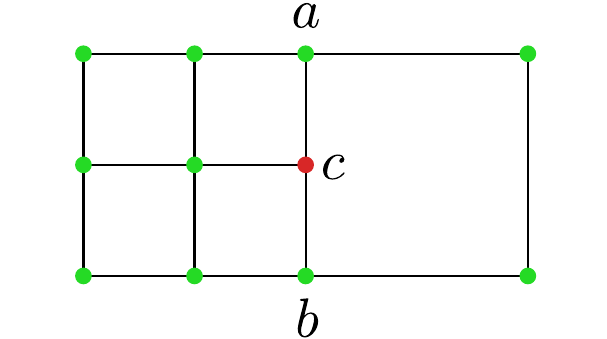}
  \hspace{1.2cm}
  \includegraphics[width=0.3\textwidth]{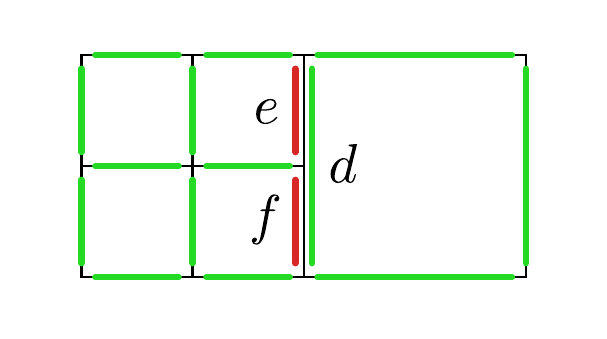}
\end{center}
\caption{Illustration of conformity constraints for lowest order elements
in 2D.  Left: nodal elements ($H^1$ subspace), constraint $c = (a+b)/2$.
Right: Nedelec elements ($\Hcurl$ subspace), constraints $e = f = d/2$.
In both cases, slave degrees of freedom are linearly interpolated from the
degrees of freedom on the master side.}
\label{fig:local-simple}
\end{figure}

To restore conformity at the coarse-fine interfaces, the degrees of freedom on
slave entities (slave DOFs) must be constrained so that the slave entities
interpolate the finite element functions of their masters. Since the masters and
slaves have matching function spaces (all elements have the same order $p$),
the slave DOFs can simply be expressed as linear combinations of the
master DOFs. Examples of such local constraints for low-order elements are shown
in Figure \ref{fig:local-simple}.

Note that in general non-conforming meshes, some DOFs on the boundaries of
master entities may themselves be constrained by another master entity, or a
chain of them. In this paper we are concerned with {\em solvable} non-conforming
meshes, where each such dependency chain ends with unconstrained degrees of
freedom. The final set of unconstrained DOFs is called \emph{true} DOFs, denoted
$u_h$. A solution vector $\widehat{u}_h$ in the partially conforming space
$\widehat{V}_h$ can thus be written as
\begin{equation}
  \widehat{u}_h = \begin{pmatrix} u_h \\ w_h \end{pmatrix} \qquad
  u_h \in \mathbb{R}^N,\; w_h \in \mathbb{R}^{(\widehat{N}-N)},
\end{equation}
where $w_h$ represents all slave DOFs and can be evaluated at any time by a
linear interpolation $w_h = W u_h$ for some interpolation matrix $W$.  We can
also write
\begin{equation}
  \widehat{u}_h = P u_h,\qquad \text{where }
  P = \begin{pmatrix} I \\ W \end{pmatrix}.
\end{equation}
$P$ will be called the \emph{conforming prolongation matrix}. The $N$ columns of
$P$ define the conforming basis functions of $V_h$, expressed as vectors in
$\widehat{V}_h$. We also use $P$ to variationally restrict the solution
$\widehat{u}_h$ to the conforming subspace $V_h$.

Since we want to upgrade an existing finite element code to support non-conforming meshes
with as few modifications as possible, we let the code assemble the stiffness
matrix $\widehat{A}$ and load vector $\widehat{b}$ in the space $\widehat{V}_h$,
as if there were no hanging nodes in the mesh. Of course, when computed in
$\widehat{V}_h$, the discrete problem \eqref{discreteprob} and the corresponding
linear system $\widehat{A}\,\widehat{u}_h = \widehat{f}$ result in a
non-conforming solution where the slave DOFs are not constrained.
However, taking $\widehat{u}_h = P u_h$ and $\widehat{v}_h = P v_h$, the variational
formulation on $V_h$ becomes
\begin{equation}
  v_h^T P^T\!\widehat{A} P u_h = v_h^T P^T \widehat{f}
  \qquad \forall v_h \in \mathbb{R}^{N},
\end{equation}
so we can solve the smaller system
\begin{equation}
  P^T\!\widehat{A} P u_h = P^T \widehat{f}, \label{eq:system}
\end{equation}
where the slave DOFs are eliminated. We then prolongate $u_h$ to obtain
a conforming solution $\widehat{u}_h \in \widehat{V}_h$.

An important point regarding essential boundary conditions in $V_h$ is that they
have to be applied in the system \eqref{eq:system} {\em after} the variational
transformation, i.e., based on the matrix $P^T\!\widehat{A}P$ and right hand
side $P^T\widehat{f}$. This is in contrast to many non-AMR finite element
settings, where one can eliminate the essential boundary conditions on element
level. Indeed, it is straightforward to verify that for a general sub-matrix
$W$, the elimination based on $\widehat{A}$ and $\widehat{f}$ followed by
variational restriction is different than performing variational restriction and
then elimination. In our experience, ensuring that the essential boundary
conditions are eliminated at the very end is the main change needed in
applications adopting the proposed AMR approach.

Another algebraic transformation that is commonly used in practice is static
condensation, see e.g. \cite{wilson74}. Since static condensation manipulates
DOFs in the interior of the element, and such DOFs are always true, AMR and
static condensation commute, i.e., static condensation followed by AMR
variational restriction is the same as AMR reduction followed by static
condensation.  The variational restriction approach can also be used to couple
AMR with hybridization, see \cite{Dobrev17}.

\section{Constructing serial $P$} \label{sec:serialP}

We now describe the algorithm to calculate the conforming prolongation matrix
in a single-processor setting. As indicated in the previous section,
factoring $P$ out of the coupled AMR matrix makes it easier to incorporate AMR
in applications. The explicit availability of $P$ can also facilitate the
development of AMR-specific linear solvers and preconditioners. Many codes apply
$P$ without assembling it \cite{WeinzierlMehl2011}, but in our case it is more convenient to build it
explicitly due to the complex dependencies that can arise in arbitrarily irregular 3D
meshes.

The $P$ matrix is constructed for a concrete instance of a partially conforming
finite element space $\widehat{V}_h$. The procedure described below constructs
both the true DOF identity $I$ and the slave interpolation $W$ submatrices of
$P$. The first step is to determine the master-slave relations of non-conforming
edges and faces. Algorithms for that are described later. Here we assume that we
are given lists of master edges and faces and their associated slaves.

\begin{figure}
\begin{center}
  \includegraphics[width=0.43\textwidth]{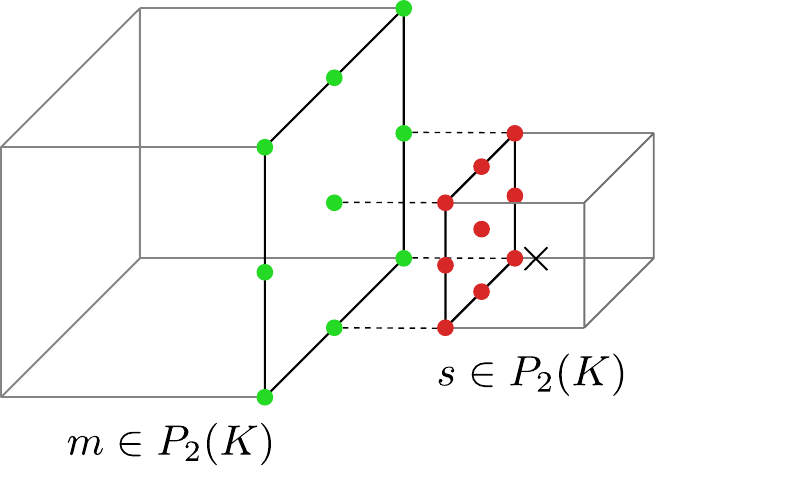}
\end{center}
\caption{Example of local constraining relations for quadratic hexahedral elements.
The slave degrees of freedom, $s$, interpolate the master degrees of freedom, $m$.
We can write $s = Q m$, where $Q$ is a local interpolation matrix
(in this case $Q \in \mathbb{R}^{9 \times 9}$).}
\label{fig:local-interp}
\end{figure}

For each slave entity $s_k$, the position of its vertices within the reference domain
of its master edge/face is known. For example, the vertices of the slave
face in Figure \ref{fig:local-interp} could have coordinates \{(0.5, 0), (1, 0),
(1, 0.5), (0.5, 0.5)\} in the reference domain of the master face. From these we
can calculate the local interpolation matrix $Q^k$ which expresses the slave
DOFs (including DOFs on the boundary of the slave face) as linear combinations
of the DOFs of the master face (including its boundary DOFs).
In the case of nodal elements, the rows of $Q^k$ contain just the shape
functions of the master face finite element evaluated in the nodal positions of
the slave face, i.e. the slave face performs nodal interpolation of the master face
finite element function. Since the polynomial degrees of the two faces match, the
interpolation is exact. Although we primarily use nodal elements, the local
interpolation matrix is not limited to nodal elements and can be computed for
any kind of finite element.

Note that sometimes a row of $Q^k$ will be of the
trivial form $Q^k_i$ = (0, ..., 0, 1, 0, ..., 0), which means that the slave DOF
is identical to the master DOF and may not be constrained at all. This is the
case, for example, for the node marked with ``$\times$'' in Figure
\ref{fig:local-interp}.
Since the matrices $Q^k$
include the boundary DOFs of the slave face, certain dependencies will be
computed multiple times --- this is done to simplify the implementation at a
relatively small computational cost.

\begin{figure}
  \begin{minipage}{.47\textwidth}
    \centering
    \includegraphics[width=0.6\textwidth]{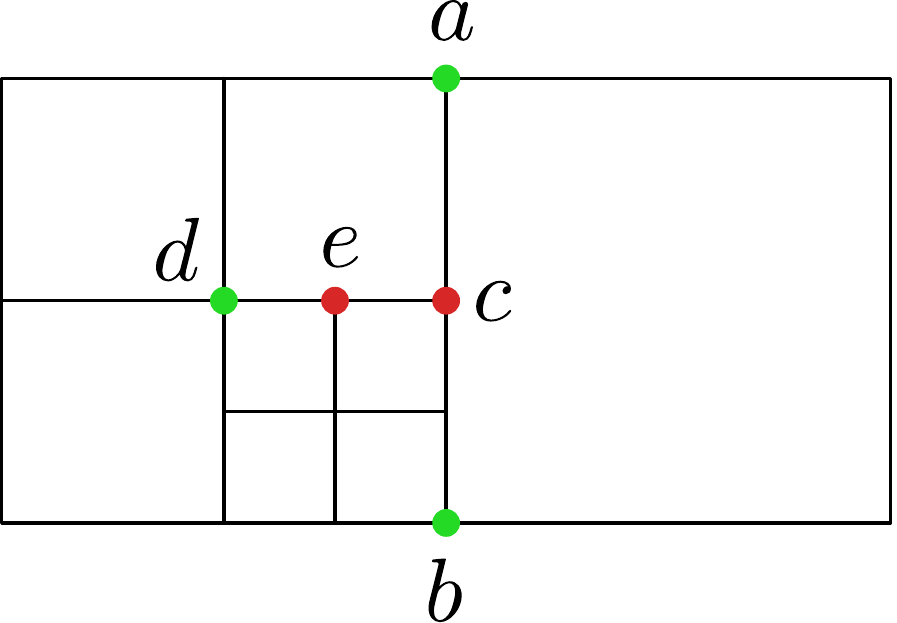}
    \caption{Example of indirect constraints: $e = (c + d)/2$ depends on
             $c = (a + b)/2$.}
    \label{fig:indirect}
  \end{minipage}
  \hfill
  \begin{minipage}{.47\textwidth}
    \centering
    \includegraphics[width=0.39\textwidth]{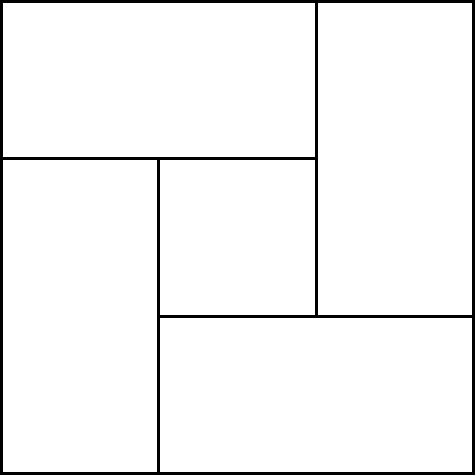}
    \caption{Unsolvable mesh with cyclic dependencies between nodal degrees of freedom.}
    \label{fig:cyclic}
  \end{minipage}
\end{figure}

As we iterate through the local interpolation matrices $Q^k$ of all slave
edges/faces, we continuously update a global dependency matrix $D_{ij}$, $1 \le i,j
\le \widehat{N}$. The matrix $D$, corresponding to the DOFs in $\widehat{V}_h$, is
represented as a sparse matrix and is initialized to identity. For each $Q^k$,
we assign the local interpolation weights to $D_{ij}$ for
each slave DOF $i$ that depends on a master DOF $j$. This allows us to identify
three kinds of rows in the dependency matrix at the end of this phase:

\medskip

\begin{enumerate}
\item {\em Identity rows}: $D_i = (0,...,0,1,0,...,0)$. These represent
  unconstrained, true degrees of freedom.

\medskip

\item {\em Direct slave rows}: non-identity $D_i = (d_{ij})$ with nonzero entries
  $d_{ij}$ corresponding only to unconstrained degrees of freedom $j$. These are
  slave degrees of freedom directly dependent on true DOFs.

\medskip

\item {\em Indirect slave rows}: the remaining non-identity rows represent slave
  degrees of freedom dependent on true and/or other slave DOFs. These rows appear
  because we do not restrict the mesh regularity, see Figure \ref{fig:indirect}
  and the discussion below.
\end{enumerate}

\medskip

Let $N$ be the number of identity rows of the $D$ matrix. This is the number of
true degrees of freedom and the dimension of the conforming finite element
space $V_h$.  We initialize $P$ as an $\widehat{N} \times N$ sparse matrix with identity
in the $N$ rows corresponding to the true DOFs.  To fill in the remaining rows,
we introduce a Boolean vector $r$ of size $\widehat{N}$, which keeps track of which rows
in $P$ have been {\em resolved}. If $i$ is a true DOF, we initialize $r_i=1$,
otherwise we set $r_i=0$.  We loop over the rows of $D$, scanning for rows where
$r_i = 0$ and $r_j \ne 0$ for all nonzero $D_{ij}$. These rows correspond to
unresolved constrained DOFs that depend on already resolved DOFs. When such
$D_i$ is found, the (sparse) row $P_i$ is set to
\begin{equation}\label{Presolve}
P_i = \sum_j D_{ij} P_j\,,
\end{equation}
and the resolution is marked by setting $r_i$ to one.

This procedure is iterated until all unresolved $P_i$ are assigned a linear
combination of already known $P_j$. In the first sweep, all direct dependencies
are resolved (e.g. the dependence of $c$ on $a$ and $b$ in Figure
\ref{fig:indirect}), which unlocks indirectly dependent slave DOFs that are
resolved in subsequent sweeps (e.g. the dependence of $e$ on $c$ in Figure
\ref{fig:indirect}). The number of sweeps needed is based on the longest
indirect dependency in the non-conforming mesh known as {\em irregularity}
\cite{Babuska80,Demkowicz89,Carstensen09}: 1-irregular meshes have only one level of
dependency, while $n$-irregular meshes will require $n$ sweeps to resolve all
rows of $P$. In practice, we typically have at most 3-irregular meshes, so only a
few sweeps are needed for all rows of $P$ to be resolved.

The algorithm is only guaranteed to finish if there are no cyclic dependencies,
i.e., if the mesh is solvable. In the serial algorithm we include an explicit
check: if there are still unresolved rows $i$ but none of them has $r_j \ne 0$
for all nonzero $D_{ij}$, we print an error message. This never happens if the
AMR process starts on an initially conforming (regular) mesh. However, it is possible to
construct a non-conforming mesh with cyclic dependencies. Figure \ref{fig:cyclic}
shows a mesh (inspired by \cite{Carstensen09}) that is representable by our data
structure but leads to an infinite loop in the unchecked $P$ matrix algorithm.
Note that in this case, a single level of additional refinement breaks the cycle
and makes the mesh solvable.

\section{Serial data structures} \label{sec:data-serial}

In this section we describe how we represent non-conforming meshes and how we
determine master-slave relations needed for the construction of $P$. A
high-order FEM code needs to be able to track vertices, edges and faces and
their relations to incident elements. Octree-based algorithms can derive this
connectivity from the structure of the tree itself \cite{WeinzierlMehl2011}, though on octree
boundaries the algorithms get complex (see interoctree connectivity in
\cite{p4est11}). Since our meshes are more general than octrees, we take a more
traditional approach with explicit links between mesh entities. A comprehensive
and mathematically elegant approach is the Sieve representation \cite{DMPlex},
storing a graph with edges connecting entities whose dimension $d$ differs by one.
However, this cell-complex representation does not generalize naturally to
non-conforming meshes \cite{DMPlexNC} and stores more information than we need.

As an alternative, we merely connect each $d>0$ dimensional entity to its
vertices and store the inverse mappings (from vertices back to edges and faces) as
hash tables. This simple data structure provides the mesh connectivity we need
and allows us to traverse the non-conforming interfaces. More specifically,

\begin{itemize}
\item edges and hanging vertices are identified by pairs of vertex indices,
\item faces are identified by four vertex indices, and contain two element references,
\item elements contain indices of up to 8 vertices, or indices of children if refined.
\end{itemize}

For a uniform hexahedral mesh, this representation requires about 290 bytes per
element, counting the vertices, edges, faces and their 32-bit indices, and
including the hash tables and the refinement trees. For a high-order code we
believe this to be acceptable, as a high-order mesh typically contains at least
8 times fewer elements than a comparable low-order mesh.

\algblockx{Struct}{EndStruct}[1]{\textbf{struct} \textsc{#1}}{\textbf{end struct}}
\algblockx{Union}{EndUnion}{\textbf{union}}{\textbf{end union}}

\begin{figure}
  \begin{minipage}{.46\textwidth}
    \centering
    \begin{algorithmic}
      \Struct{Element}
        \State \textbf{bool} refined
        \Union
          \State \textbf{int} vertex[8]
          \State \textbf{int} child[8]
        \EndUnion
        \State \textbf{int} parent
        \State \textbf{int} rank \Comment{see Section \ref{sec:parallel}}
      \EndStruct
    \end{algorithmic}
    \caption{Data structure for both refined and leaf elements (pseudocode).}
    \label{fig:elem-struct}
  \end{minipage}
  \hfill
  \begin{minipage}{.46\textwidth}
    \centering
    \includegraphics[width=0.56\textwidth]{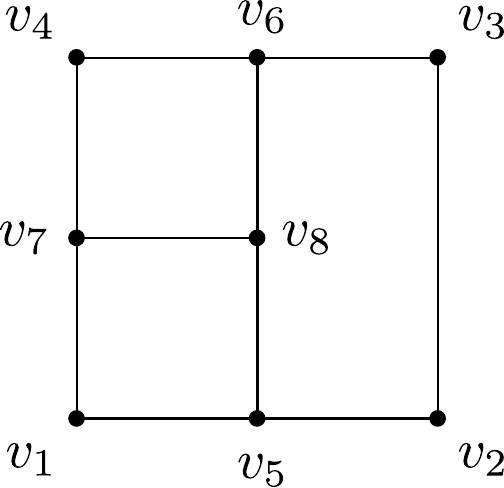}
    \caption{Example of a hexahedron face $(v1, v2, v3, v4)$, adjacent to three
             neighboring hexahedra in a non-conforming mesh.}
    \label{fig:meshstruct}
  \end{minipage}
\end{figure}

\subsection{Refinement trees} \label{sec:tree}

Since we aim for dynamic AMR capable of coarsening, we need to store the
refinement history. When the elements of the initial mesh get refined, they are
removed from the list of active elements and serve as roots of refinement trees.
A refined element points to its children instead of its vertices, as shown in
pseudocode in Figure \ref{fig:elem-struct}. To save memory, we use 32-bit
integers to identify entities instead of pointers. An element can have fewer
than 8 children if it is refined anisotropically. We store a computed list
of all leaf nodes of all refinement trees. This list represents the current
state of the mesh.

\subsection{Hash maps} \label{sec:maps}

The mapping from vertices to other entities is done by hash tables.
We illustrate this on an example shown in Figure \ref{fig:meshstruct}.
Assume $(v_1, v_2, v_3, v_4)$ is a 3D face. Across this face there are three
neighboring elements with faces $(v_1, v_5, v_8, v_7)$, $(v_5, v_2, v_3, v_6)$
and $(v_7, v_8, v_6, v_4)$. To access edges and faces from the incident elements,
we can call

\begin{algorithmic}
   \State $e_1 \gets$ GetEdge$(v_1, v_5)$
   \State $e_2 \gets$ GetEdge$(v_5, v_2)$
   \State $f_1 \gets$ GetFace$(v_1, v_5, v_8, v_7)$
   \State $f_2 \gets$ GetFace$(v_5, v_2, v_3, v_6)$
   \State ...
\end{algorithmic}

A similar hash map exists for hanging vertices and is used
to traverse non-conforming interfaces in Section \ref{sec:masterslave}:

\begin{algorithmic}
   \State $v_5 \gets$ GetVertex$(v_1, v_2)$
   \State $v_6 \gets$ GetVertex$(v_4, v_3)$
   \State $v_7 \gets$ GetVertex$(v_1, v_4)$
   \State $v_8 \gets$ GetVertex$(v_5, v_6)$
\end{algorithmic}

Each Get function searches the respective hash table and returns a 32-bit
entity index. The order of arguments is not relevant, i.e., GetVertex$(v_1,
v_2)$ = GetVertex$(v_2, v_1)$.
We use the convention that the Get functions create the requested entity if it
does not exist. To query its existence, each GetX function has a FindX
counterpart taking the same arguments and returning {\bf nil} (encoded as $-1$)
if the entity does not exist, instead of creating it. We make sure the hash tables
are properly (dynamically) sized to get amortized $O(1)$ complexity of each query.

\subsection{Refinement and coarsening} \label{sec:refinement}

To refine e.g.\ a hexahedron, we first create new mid-edge vertices (or get already
existing ones) by calling GetVertex$(\cdot, \cdot)$ for all 12 edges. Six
mid-face vertices are accessed through opposite mid-edge vertices; the two
options per face must be checked first with FindVertex$(\cdot, \cdot)$ to
prevent creating a duplicate mid-face vertex, in case it already exists. The
mid-element vertex is accessed from any two opposite mid-face vertices. New
elements are created with indices of the new vertices. The old element is
marked as refined and its vertex indices are replaced with indices of the child
elements.
Reference counting is used to keep track of the number of elements pointing to
each vertex, edge and face. If the reference count drops to zero, the entity is
deleted.

\begin{figure}
\begin{center}
  \includegraphics[width=0.85\textwidth]{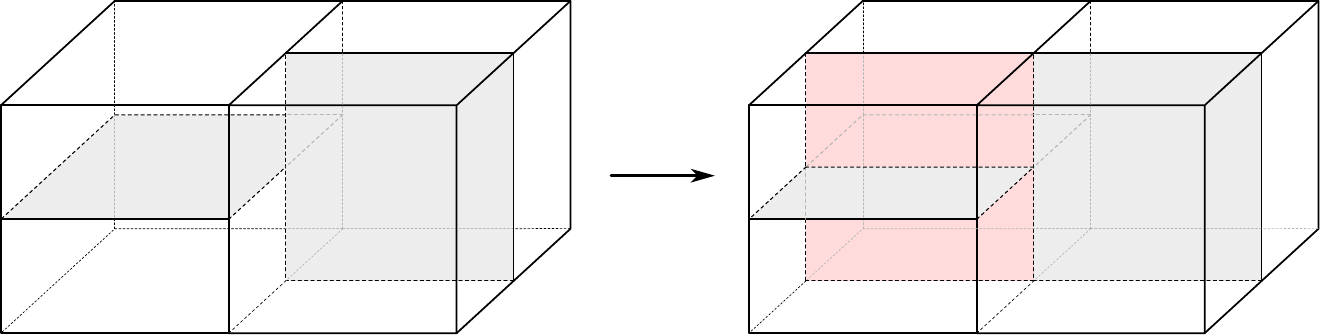}
\end{center}
\caption{Anisotropic refinements may lead to an inconsistent mesh (left). Additional
  refinements are needed to restore valid master-slave relations between faces
  (right). }
\label{fig:aniso}
\end{figure}

Our data structure allows for anisotropic refinement. The challenge is
maintaining consistency of the mesh when refining anisotropically. It is easy to
refine adjacent elements in such a way that their faces are not subsets or
supersets of each other, so there is no way to determine a master-slave
relationship. This is illustrated in Figure \ref{fig:aniso}. To resolve such a
situation, the neighbor of the element being refined has to be refined as well.
Such {\em forced} refinements may propagate, until a globally consistent mesh is
obtained. We defer the treatment of this ripple effect to a follow-up
publication.

When coarsening an element, the indices of its original vertices are first
retrieved from the children (since the refinement pattern is known). All of the
children need to be leaves, otherwise they are recursively coarsened first. In
an inverse process to refinement, the child elements are destroyed and
associated vertex, edge and face reference counts are decremented. The coarse
element is reactivated, and the appropriate reference counts are incremented.

\subsection{Calculating master-slave relations} \label{sec:masterslave}

We loop over all leaf elements and their faces. Each face, triangular
or quadrilateral, is checked if it is a master face by invoking the function
TraverseTriFace (Algorithm \ref{alg:traverse-tri-face}) or TraverseQuadFace
(Algorithm \ref{alg:traverse-quad-face}). These functions attempt to recursively
descend an implicit face refinement ``tree'' and return a list of slave faces,
if they can be reached, or an empty list if the starting face is conforming.
The function FindVertex is used here to obtain vertices that were created at the
center of edges of previously refined elements. During descent, the position
within the master face's reference domain is tracked by the points $p_i$, which
are initialized to the reference domain corners (at $depth = 0$).

In the case of quadrilateral faces, the face refinement ``tree'' is traversed as
a binary tree, i.e., isotropically refined faces are treated as if refined
anisotropically twice (horizontally and then vertically relative to their
reference domain, or vice versa). The function FaceSplitType (Algorithm
\ref{alg:face-split-type}) determines which way a face is split, again using
FindVertex to look up mid-edge and mid-face vertices (refer to Figure
\ref{fig:meshstruct} for a quadrilateral face example).

\begin{algorithm}
  \caption{List slave faces of a potential triangular master face.}
  \label{alg:traverse-tri-face}
  \begin{algorithmic}[1]
    \Function{TraverseTriFace}{$v_1, v_2, v_3,\, p_1, p_2, p_3,\, depth$}
      \State $v_{12} \gets$ FindVertex($v_1, v_2),\ p_{12} \gets (p_1 + p_2)/2$
      \State $v_{23} \gets$ FindVertex($v_2, v_3),\ p_{23} \gets (p_2 + p_3)/2$
      \State $v_{31} \gets$ FindVertex($v_3, v_1),\ p_{31} \gets (p_3 + p_1)/2$
      \State $list \gets [\,]$
      \If{($v_{12}$ \textbf{not nil}) \textbf{and}
          ($v_{23}$ \textbf{not nil}) \textbf{and}
          ($v_{31}$ \textbf{not nil})}
        \State $list \gets \mbox{TraverseTriFace}(v_1, v_{12}, v_{31},
               p_1, p_{12}, p_{31},\, depth+1)\ +$
        \State \hspace{10.8mm}$\mbox{TraverseTriFace}(v_{12}, v_2, v_{23},\,
               p_{12}, p_2, p_{23},\, depth+1)\ +$
        \State \hspace{10.8mm}$\mbox{TraverseTriFace}(v_{31}, v_{23}, v_3,\,
               p_{31}, p_{23}, p_3,\, depth+1)\ +$
        \State \hspace{10.8mm}$\mbox{TraverseTriFace}(v_{12}, v_{23}, v_{31},\,
               p_{12}, p_{23}, p_{31},\, depth+1)$
      \ElsIf{$depth > 0$}
        \State $slave \gets$ FindFace($v_1, v_2, v_3$, \textbf{nil})
        \If{$slave$ \textbf{not nil}}
          \State $list \gets [\{slave, p_1, p_2, p_3\}]$
        \EndIf
      \EndIf
      \State \Return $list$
    \EndFunction
  \end{algorithmic}
\end{algorithm}

\begin{algorithm}
  \caption{List slave faces of a potential quadrilateral master face.}
  \label{alg:traverse-quad-face}
  \begin{algorithmic}[1]
    \Function{TraverseQuadFace}{$v_1, v_2, v_3, v_4, p_1, p_2, p_3, p_4, depth$}
      \State $split \gets$ FaceSplitType($v_1, v_2, v_3, v_4$)
      \State $list \gets [\,]$
      \If{$split$ = Vertical}
        \State $v_{12} \gets$ FindVertex($v_1, v_2),\ p_{12} \gets (p_1 + p_2)/2$
        \State $v_{34} \gets$ FindVertex($v_3, v_4),\ p_{34} \gets (p_3 + p_4)/2$
        \State $list \gets \mbox{TraverseQuadFace}(v_1, v_{12}, v_{34}, v_4,
               p_1, p_{12}, p_{34}, p_4, depth+1)\ +$
        \State \hspace{10.8mm}$\mbox{TraverseQuadFace}(v_{12}, v_2, v_3, v_{34},
               p_{12}, p_2, p_3, p_{34}, depth+1)$
      \ElsIf{$split$ = Horizontal}
        \State $v_{23} \gets$ FindVertex($v_2, v_3),\ p_{23} \gets (p_2 + p_3)/2$
        \State $v_{41} \gets$ FindVertex($v_4, v_1),\ p_{41} \gets (p_4 + p_1)/2$
        \State $list \gets \mbox{TraverseQuadFace}(v_1, v_2, v_{23}, v_{41}, p_1,
               p_2, p_{23}, p_{41}, depth+1)\ +$
        \State \hspace{10.8mm}$\mbox{TraverseQuadFace}(v_{41}, v_{23}, v_3, v_4,
               p_{41}, p_{23}, p_3, p_4, depth+1)$
      \ElsIf{$depth > 0$}
        \State $slave \gets$ FindFace($v_1, v_2, v_3, v_4$)
        \State $list \gets [\{slave, p_1, p_2, p_3, p_4\}]$
      \EndIf
      \State \Return $list$
    \EndFunction
  \end{algorithmic}
\end{algorithm}

\begin{algorithm}
  \caption{Determine whether a face is split horizontally, vertically or not split.}
  \label{alg:face-split-type}
  \begin{algorithmic}[1]
    \Function{FaceSplitType}{$v_1, v_2, v_3, v_4$}
      \State $v_{12} \gets$ FindVertex$(v_1, v_2)$
      \State $v_{23} \gets$ FindVertex$(v_2, v_3)$
      \State $v_{34} \gets$ FindVertex$(v_3, v_4)$
      \State $v_{41} \gets$ FindVertex$(v_4, v_1)$
      \State $midf_1 \gets$ ($v_{12}$ \textbf{not nil and} $v_{34}$
             \textbf{not nil}) ? FindVertex($v_{12}, v_{34}$) : \textbf{nil}
      \State $midf_2 \gets$ ($v_{23}$ \textbf{not nil and} $v_{41}$
             \textbf{not nil}) ? FindVertex($v_{23}, v_{41}$) : \textbf{nil}
      \If{$midf_1$ \textbf{is nil and} $midf_2$ \textbf{is nil}}
        \State \Return NotSplit
      \Else
        \State \Return ($midf_1$ \textbf{not nil}) ? Vertical : Horizontal
      \EndIf
    \EndFunction
  \end{algorithmic}
\end{algorithm}

An analogous function, TraverseEdge, is employed to construct the list of master-slave
relations for edges. Edges of each leaf element are checked for slave edges by
traversing the local vertex hierarchy in a similar recursive manner.

The 2D case is treated with the same code as a degenerate case: triangular and
quadrilateral elements merely have a different number of vertices/edges in
internal geometry tables and they have zero faces, so creation and traversal of
faces is skipped.

We make no assumptions regarding the level of refinement of adjacent elements.
The simple algorithms above work for general non-conforming 3D meshes, including
meshes with aniso\-tro\-pic refinements.

\subsection{Element neighbors} \label{sec:neighbors}

We define two elements to be neighbors if their closed sets have a
non-empty intersection, even if the intersection is just a single vertex (i.e.
we consider all of vertex-, edge- and face-neighbors). Although algorithms
exist to determine such neighbors from the structure of the refinement tree
itself, the handling of multiple trees may become complex \cite{p4est11}.

Instead, we employ an approach based on a Boolean matrix of element-vertex
connectivity. Let $B$ be a Boolean matrix containing one row per element and one
column per mesh vertex. For each element, its row marks vertices that belong to
the element, including possible mid-edge and mid-face vertices of neighboring
refined elements. The incident vertices are collected for each element by a
procedure similar to Algorithms \ref{alg:traverse-tri-face} and
\ref{alg:traverse-quad-face}. Then, if $e$ is a Boolean vector representing the
set of elements whose neighbors we need to find, the vector $e' = BB^Te$ will
represent the original set of elements plus all its neighbors. The $BB^T$
product is not stored, only its action is computed on demand. To save even more
memory, the obvious corner vertices of each element are also not explicitly
stored in $B$. The matrix is thus empty if the mesh is conforming.

\section{Parallelization} \label{sec:parallel}

We split the mesh into $K$ disjoint regions (see next section), where $K$ is the
number of MPI tasks. The partitioning is element based, i.e., each element is
assigned to one of the tasks $k$. The vertices, edges and faces on the boundary
of each task's region are duplicated, so that each region can be treated as an
isolated mesh and passed to the serial part of the code. Repeating the approach from Section
\ref{sec:amr}, we assemble the stiffness matrices $\widehat{A}_k$ and load
vectors $\widehat{f}_k$ locally on each MPI task $k$, as if the mesh was conforming
and serial on each task.

Globally, we have a new, parallel finite element space $\widetilde{V}_h \supset
\widehat{V}_h$, with $\widetilde{N} = \dim{\widetilde{V}_h} >
\dim{\widehat{V}_h} = \widehat{N}$. The parallel stiffness matrix
$\widetilde{A}$ consists of diagonal blocks $\widehat{A}_k$ and the parallel
load vector $\widetilde{f}$ contains the blocks $\widehat{f}_k$. Again, the
solution of the parallel system $\widetilde{A}\widetilde{u}_h = \widetilde{f}$
would be disconnected along the interfaces of the $K$ regions, but we can use
the same variational restriction approach to obtain a globally conforming
parallel solution. The parallel conforming prolongation matrix, still denoted
$P$, is now of the size $\widetilde{N} \times N$ and conveniently handles both
conforming interpolation and parallel decomposition, by mapping from the space
$V_h$ directly to $\widetilde{V}_h$.

In our current implementation, we explicitly form the parallel triple matrix
product $P^T\!\widetilde{A}P$ using the $RAP$ triple-matrix-product kernel
from the {\em hypre} library \cite{hypre}.  This kernel is highly optimized in
{\em hypre} as it is used internally for the coarse-grid operator construction
in its algebraic multigrid solvers. For higher order elements we are also
considering direct evaluation of the action of $P^T\!\widetilde{A}P$, without
assembly.

\subsection{Tree-based partitioning} \label{sec:partitioning}

Similarly to \cite{GriebelZumbusch1997,Campbell2003,TOG2005,Mitchell07,p4est11,Zoltan12}, we partition the
mesh by splitting a space-filling curve (SFC) obtained by enumerating
depth-first all leaf elements of all refinement trees. We assume that the
elements of the coarse mesh are ordered as a sequence of face-neighbors, so that
a globally continuous SFC can be obtained. For ordering the coarse mesh (i.e.,
the refinement tree roots) we use the Gecko library \cite{Gecko}. Once the
coarse mesh is ordered we invoke depth-first traversal for the ordered and
oriented roots. Orienting a root element means assigning an initial traversal
state (see below) so that SFCs of successive trees are connected.

\begin{figure}
  \begin{minipage}{.45\textwidth}
    \centering
    \begin{tabular}{@{}ccc@{}}
      \toprule
      state & child order & child state\\
      \midrule
      0 & (0, 1, 2, 3) & (1, 0, 0, 5) \\
      1 & (0, 3, 2, 1) & (0, 1, 1, 4) \\
      2 & (1, 2, 3, 0) & (3, 2, 2, 7) \\
      3 & (1, 0, 3, 2) & (2, 3, 3, 6) \\
      4 & (2, 3, 0, 1) & (5, 4, 4, 1) \\
      5 & (2, 1, 0, 3) & (4, 5, 5, 0) \\
      6 & (3, 0, 1, 2) & (7, 6, 6, 3) \\
      7 & (3, 2, 1, 0) & (6, 7, 7, 2) \\
      \bottomrule
    \end{tabular}
    \caption{Ordering of subtrees to obtain the Hilbert curve in 2D.}
    \label{tab:hilbert}
  \end{minipage}
  \hfill
  \begin{minipage}{.44\textwidth}
    \centering
    \includegraphics[width=0.45\textwidth]{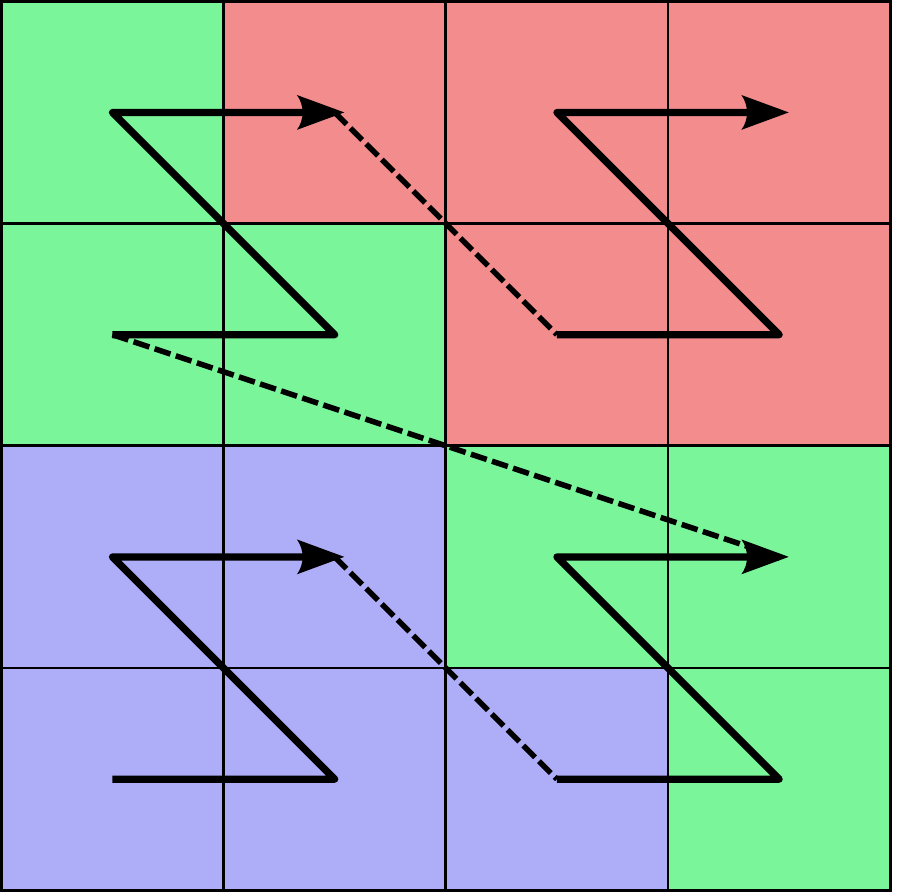}
    \hfill
    \includegraphics[width=0.45\textwidth]{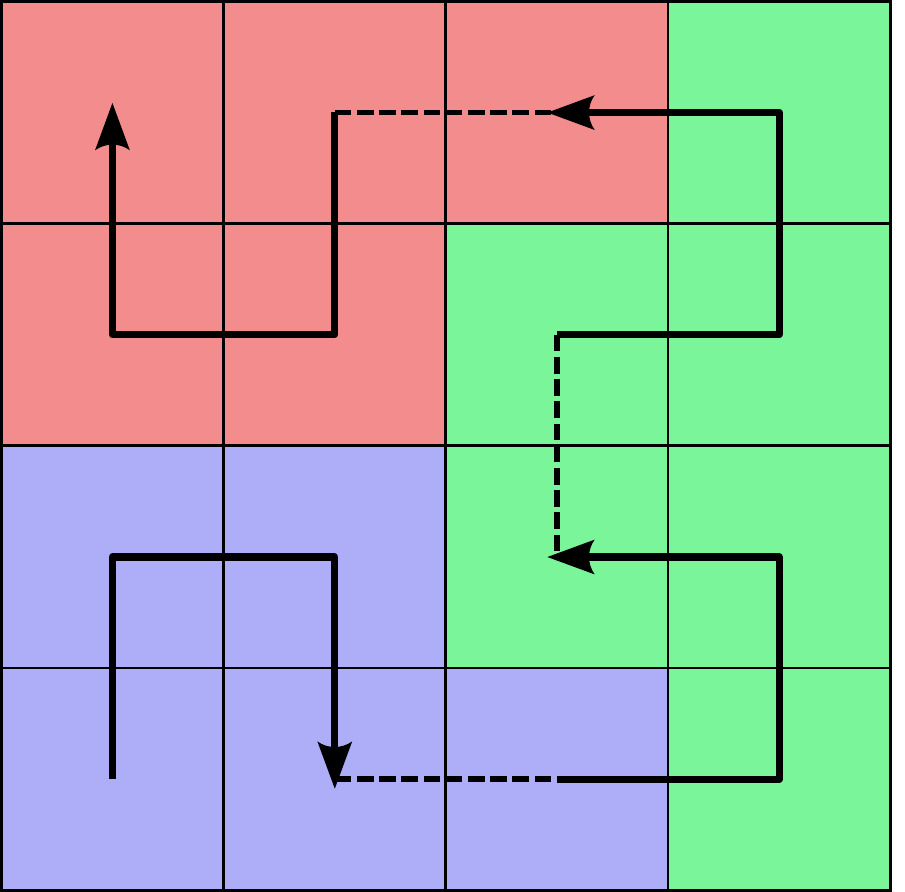}
    \caption{Example of Z-curve (left) vs. Hilbert curve partitioning.}
    \label{fig:partition}
  \end{minipage}
\end{figure}

In the case of quadrilaterals and hexahedra, the simplest traversal with a fixed
order of children at each tree node leads to the well-known Z-curve, as shown in
Figure \ref {fig:partition} (left). However, it turns out that simply by
changing the order the subtrees are visited at each internal node, the Hilbert
curve can be obtained at no additional cost (i.e., the elements do not need to
be physically rotated). The Hilbert curve produces continuous partitions
and leads to better interprocess connectivity \cite{Campbell2003} than the
Z-curve. In the simpler 2D case, each tree node can be in 8 states (4
rotations times 2 inversions). Given a state of the root node, we descend to
subtrees in the appropriate order, passing a new state to each, according to
Table \ref{tab:hilbert}. Example of a Hilbert curve partitioning is shown in
Figure \ref{fig:partition} (right). The method is analogous for hexahedra, where
a tree node has 24 states.

Generating a continuous SFC is more difficult for triangles, prisms and in the
presence of anisotropic refinements. The Sierpinski curve is optimal for
triangles \cite{Mitchell07}, but we would have to switch to triangle bisection.
Prisms admit a continuous SFC analogous to hexahedra, but not all elements can
be face neighbors (only edge neighbors). The most serious difficulty is with
anisotropic refinement, which probably prevents a continuous SFC, and generating
a good leaf sequence is an open problem.

At the beginning of the computation, we assume that the mesh is small enough to
fit completely in the memory of each MPI task. Each task uses a copy of the
serial mesh and traverses its trees as described. The list of leaf elements is
split into equal sized parts and each part is assigned to one MPI task.
Each task thus owns a region of the mesh. The region may not have a minimal
surface, but it should be continuous and relatively compact due to the SFC
locality \cite{Campbell2003}.

\subsection{Ghost layer} \label{sec:ghost-layer}

Although the serial code sees a non-overlapping decomposition of the mesh, we
internally hold the (pruned) refinement trees on all MPI tasks, as in \cite{p4est11}.
For the purpose of constructing $P$, we also track \emph{ghost
elements} for each MPI task. A ghost element is one that is a vertex-, edge-, or
face-neighbor to an element owned by an MPI task, but is itself owned by another
task. The set of all ghost elements forms a minimal layer of elements enclosing
a task's region. This layer needs to be kept synchronized with neighboring tasks.

\begin{figure}
\begin{center}
  \includegraphics[width=0.99\textwidth]{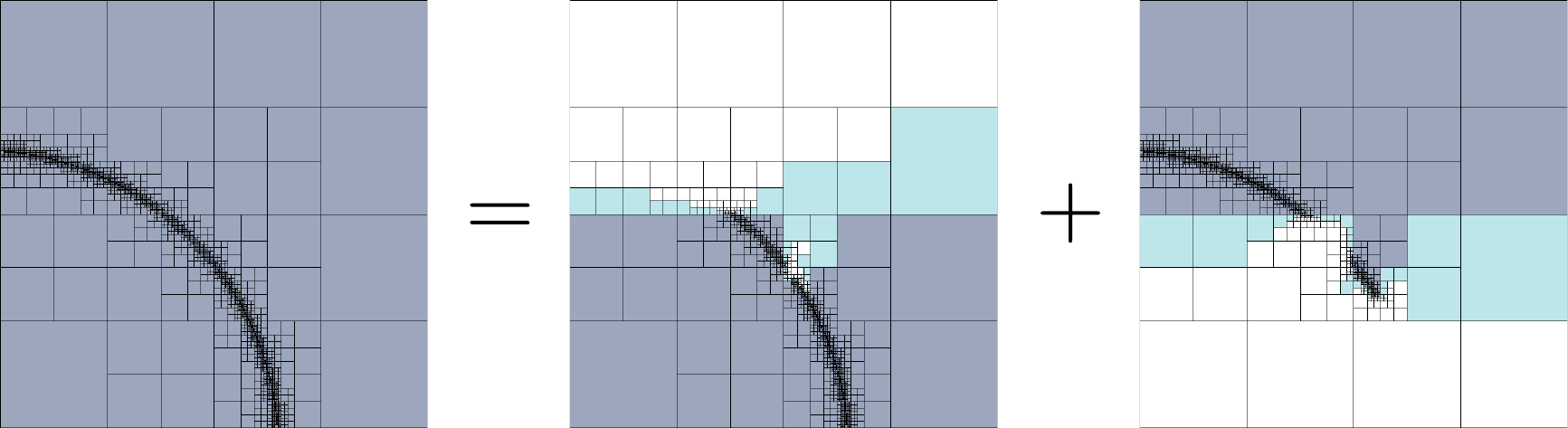}
\end{center}
\caption{A 2D mesh partitioned by the Hilbert curve between two MPI tasks.
Light-gray elements are ghost elements. White areas have been pruned and do not
represent real elements.}
\label{fig:ghosts}
\end{figure}

On each task,
elements beyond the ghost layer may not correspond to real elements, and we prune
the refinements in this area immediately after partitioning the serial mesh and
identifying the ghost elements. Pruning the refinement tree means removing
(coarsening) all subtrees that contain only leaves not owned by the current task
nor belonging to the ghost layer \cite{BBH2011}. After this step, the mesh becomes
fully distributed as no single task sees all of the actual leaf elements. This is
illustrated in Figure \ref{fig:ghosts} for two tasks in 2D.

The last type of elements we identify are parallel boundary elements. This layer
is a subset of the region owned by each MPI task formed by elements which are
vertex-, edge-, or face-neighbors to the ghost layer. The boundary elements
constitute ghost elements of one or more neighboring tasks and any change in
them needs to be communicated to the neighbors to keep all the ghost layers
synchronized. For example, if a refinement occurs in one task's boundary layer,
a message needs to be sent to all neighboring tasks that share the elements, so
they can duplicate the refinement in their ghost layer.

\subsection{Load balancing} \label{sec:balancing}

Balancing the mesh so that each MPI task has the same number of elements ($\pm1$
if the total number of elements is not divisible by the number of tasks) is
relatively straightforward in the context of SFC based partitioning. A~parallel
partial sum (\code{MPI\_Scan}) is done to determine new beginnings of the
partitions in the global sequence of leaf elements. New assignments of
\code{Element::rank} within each parallel region can then be computed. Next, the
new assignments within the ghost layers must be communicated between neighboring
tasks (this is the most expensive part). This facilitates the final step, in
which elements no longer assigned to the current task are sent to the new
owners, together with a layer of ghost elements, to ensure a valid new state of
the distributed mesh. The exchange of elements ends when each task obtains the
right number of elements, which is known beforehand. After the load balancing
procedure, a pruning step is done on each task to remove branches of the
refinement tree that no longer need to be represented, i.e., subtrees that only
contain leafs beyond the ghost layer are removed.

\subsection{Message encoding} \label{sec:encoding}

The prime advantage of storing the coarse mesh (i.e., all refinement tree roots)
on all MPI tasks is the ability to refer to any element of the global mesh by
the index of a root element and the refinement path to the element in question.
However, rather than identify elements in this way individually, we develop an
algorithm to encode a {\it subtree} of the global refinement hierarchy, since
our messages usually carry information for multiple elements at once.

Given a set of leaf elements in a local mesh, we use the member \code{
Element::parent} recursively to identify the set of trees participating in the
element set and store their indices. In each tree, we then descend along the
same paths back to the leaves. At each node, we output an 8-bit mask indicating
which subtrees contain elements from the set. We output a zero mask to terminate
the descent at leaves.

The functionality of encoding element sets is used for load balancing, but also
as a~basis for encoding sets of vertices, edges, faces and their degrees of
freedom, as required in Section \ref{sec:parallelP}. A~vertex, edge or face can
be identified by the index of an element in an element set, followed by the
local number of the vertex/edge/face within that element. A DOF is
identified by its index within its mesh entity.

We use variable length \code{MPI\_BYTE} messages and try to collect as much data
as possible for communication between each pair of tasks in order to minimize
the total number of messages.

\subsection{Constructing parallel $P$} \label{sec:parallelP}

The algorithm to build the $P$ matrix in parallel is more complex, but
conceptually similar to the serial algorithm. We still express slave DOF rows of
$P$ as linear combinations of other rows; however, some of them may be located on
other MPI tasks and must be communicated first.

The ghost layer allows us to determine, without
communication, which vertices, edges and faces are shared by more than one
task. As we traverse the local mesh and its ghost layer using the serial
algorithms of Section \ref{sec:masterslave}, we build a list of groups of MPI
ranks, which we call communication groups. Each vertex, edge and face is
assigned a communication group index, and an owner rank. In a group, the task
with the lowest rank is defined to be the owner. We store the owner explicitly
for each mesh object, because as the next step, we modify the communication
groups so that master and slave faces and edges are grouped together. This is to
ensure that information will flow from the owners and masters to non-owners and
slaves.

Next, for a concrete finite element space, we assign $\widehat{N}_k$ local
degrees of freedom on each MPI task $k$. These are the DOFs of all elements of
the task's region, whether slave, master or conforming. Globally, $\sum
\widehat{N}_k = \widetilde{N}$. Each local DOF has a unique owner, inherited
from its underlying vertex, edge or face. DOFs with a remote owner are expected
to later receive a row of the $P$ matrix that will make them identical to a
remote true DOF.

In addition to local DOFs, we also assign $\widehat{G}_k$ virtual DOFs within
the ghost layer (if not already marked as local), all with indices greater
than $\widehat{N}_k$. We define vectors $o$, $g$, $r$ of size $\widehat{N}_k +
\widehat{G}_k$, where for each DOF $i$: $o_i$ is the owner rank, $g_i$ is the
communication group index and $r_i$ is the resolution flag, as in Section
\ref{sec:serialP}.

On each task we define a local dependency matrix $D_{ij}$, this time $1 \le i
\le \widehat{N}_k$ and $1 \le j \le \widehat{N}_k + \widehat{G}_k$, initialized to identity on the $\widehat{N}_k$
block. Reusing much of the serial code, we again collect the slave interpolation
weights $D_{ij}$ for each slave DOF $i$ that depends on master DOFs $j$. Note
that $j$ may now be a ghost DOF, $j > \widehat{N}_k$. Still without any communication, we
can identify the following types of rows $D_i$:

\medskip

\begin{enumerate}
\item {\em Identity rows}, with $o_i = k$. These are true DOFs that we own.

\medskip

\item {\em Identity rows}, with $o_i \ne k$. Here $i$ is a true DOF owned by
  another rank.

\medskip

\item {\em Local slave rows}, non-identity $D_i$ with $j \le \widehat{N}_k$ for
  all nonzero $D_{ij}$.  These are slave DOFs dependent only on local DOFs.

\medskip

\item {\em Remote slave rows}, non-identity $D_i$ with at least one nonzero
  $D_{ij}$ for $j > \widehat{N}_k$. To resolve these slave DOFs, one or more remote $P$
  rows need to be received first.
\end{enumerate}

\medskip

\begin{algorithm}[t]
  \caption{Construct parallel $P$ matrix.}
  \label{alg:parallelP}
  \begin{algorithmic}[1]
    \Function{ConstructParallelP}{mesh, space}
      \State assign $\widehat{N}_k + \widehat{G}_k$ local DOFs
      \State calculate local dependency matrix $D_{ij}$ and vectors $o_i$, $g_i$
      \State $N_k \gets$ number of true DOFs owned
      \State $N_{0,k} \gets$ \code{MPI\_Scan} on $N_k$ \Comment{{\it first column
             of our $P$ partition}}
      \State $r_i \gets 0$ for all $i$
      \ForAll{$N_k$ owned true DOFs $i$}
        \State $P_i \gets $ identity, $1$ at global column $N_{0,k} + j$,
               $\ 0 \le j < N_k$
        \State $r_i \gets 1$
        \State add $P_i$ to outgoing messages to ranks $r \ne k$ of group $g_i$
      \EndFor
      \Repeat \Comment{{\it main outer loop}}
        \State \code{MPI\_Isend} all nonempty outgoing messages
        \While{\code{MPI\_Iprobe} for some rank $r$}
          \State decode message from $r$, in particular the DOF numbers $i$
          \ForAll{decoded DOFs $i$}
            \State $P_i \gets$  row from message
            \State $r_i \gets 1$
            \State add $P_i$ to outgoing messages to other ranks if needed (see
                   text) \label{alg:parallelP:forward}
          \EndFor
        \EndWhile
        \While{exists $i$ s.t. $r_i = 0$ and $r_j \ne 0$ for nonzero $D_{ij}$}
               \Comment{{\it inner loop}}
          \State $P_i \gets \sum_j D_{ij} P_j$
          \State $r_i \gets 1$
          \State add $P_i$ to outgoing messages to ranks $r \ne k$ of group $g_i$
        \EndWhile
      \Until{$r_i$ = 1 for all $i$}
      \State \Return $P$
    \EndFunction
  \end{algorithmic}
\end{algorithm}

Let $N_k$ be the number of rows of type 1 on task $k$. We perform a parallel partial
sum (\code{MPI\_Scan}) on $N_k$ to determine the global index $N_{0,k}$ of the first true
DOF for each rank, i.e. the column partitioning of $P$. We then run a version of
the iterative algorithm of Section \ref{sec:serialP}, wrapped in one more outer
loop. The outer loop communicates remote rows of $P$ when there is nothing more
to be done by the {\em serial} (local) inner loop. In the first outer iteration,
the inner loop resolves all DOFs that do not depend on remote DOFs. Then we send
the rows of all DOFs that have just been resolved to neighboring tasks. The
group $g_i$ of a DOF determines which ranks need to be sent its $P$ row. We
collect the rows and send a single message per neighbor. The messages are
decoded by the recipients and the received rows, containing global column
numbers, are assigned to appropriate ghost DOFs. This unlocks (by setting the
ghost DOF $r_i$) in the next outer iteration more DOFs that have been waiting
for remote data. The complete procedure is summarized in Algorithm
\ref{alg:parallelP}.

\begin{figure}
  \centering
  \includegraphics[width=0.35\textwidth]{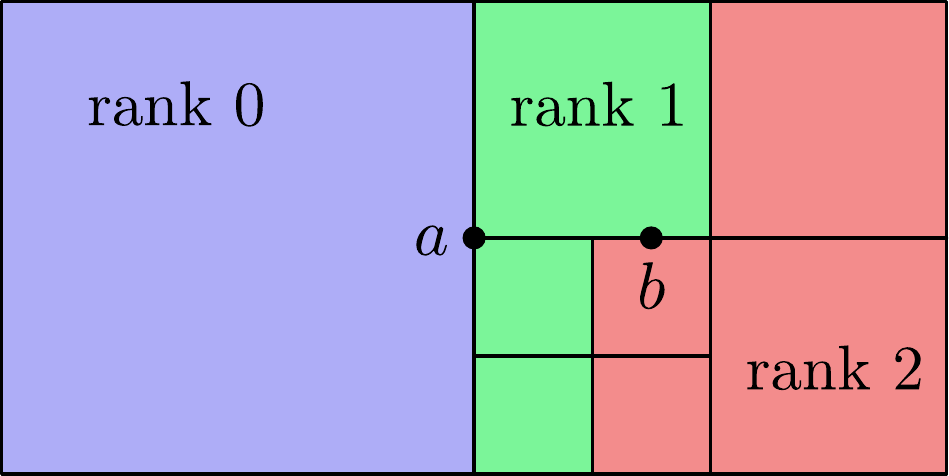}
  \caption{Message from rank~0 to rank~1 needs to be forwarded to rank~2.}
  \label{fig:forwarding}
\end{figure}

Situations may occur in which a ghost layer in an MPI task does not provide
complete information on the ranks that depend on a particular master DOF. One
such setup is shown in Figure~\ref{fig:forwarding}, where DOF $a$ is owned by
rank~0, and the corresponding row $P_a$ is needed by both rank~1 and rank~2
(e.g., to constrain DOF $b$ owned by rank 2). The group $g_a$ should therefore
be $\{0, 1, 2\}$, but rank~0 does not see rank~2 elements---these are beyond its
ghost layer. Since rank 0 sees the group $\{0,1\}$, the message is only sent to
rank~1. However, rank~1 sees the correct group, and so it forwards the message
to rank~2. This is accomplished by line~\ref{alg:parallelP:forward} of
Algorithm~\ref{alg:parallelP}, which compares the sender's version of the group
(encoded in the message) with the recipient's own view, and the message is sent
to the missing ranks. To prevent infinite loops, a message is only allowed to be
forwarded once. Testing shows that the number of rows that need forwarding is
rarely more than 1\%, typically in cases where the number of elements per task
is small.

\section{Numerical results} \label{sec:numerical}

In this section we present numerical experiments that illustrate the performance
of our unstructured AMR algorithms in practice. The results were obtained with
the open-source implementation of the proposed methods in the MFEM finite element library
\cite{mfem}.

\subsection{Model problems} \label{sec:model}

As a first test, we solve a~standard AMR benchmark problem with a known exact
solution, described in \cite{Mitchell2011}. The goal is to reveal potential
errors in the prolongation matrix algorithm by checking the behavior of the
convergence curves. We also examine the benefit of AMR compared to uniform
refinement.

The problem is the Poisson's equation $-\Delta u = f$ on the unit square $\Omega
= (0,1)^2$ with a Dirichlet boundary condition. The right hand side $f$ is
chosen so that the exact solution is
\begin{equation}
u(x, y) = \arctan\left(\alpha\left(
   \sqrt{(x - x_c)^2 + (y - y_c)^2} - r \right)\right).
\end{equation}
The
solution has a sharp circular wave front of radius $r$ centered at $(x_c, y_c)$,
as shown in Figure \ref{fig:layer2d} (left). The parameters are $\alpha = 200$,
$(x_c, y_c) = (-0.05, -0.05)$ and $r = 0.7$.  We start with a mesh consisting of
$4\times 4$ squares and use nodal finite elements of degree $p = 1,2,4,8$ (a
subspace of $H^1$) to discretize the problem. We then perform an adaptation loop
consisting of the following steps:

\begin{figure}
  \centering
  \includegraphics[width=0.3\textwidth]{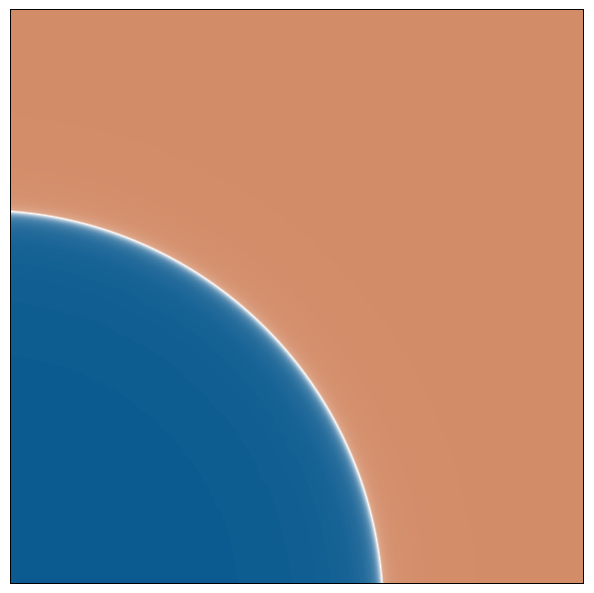} \hfill
  \includegraphics[width=0.3\textwidth]{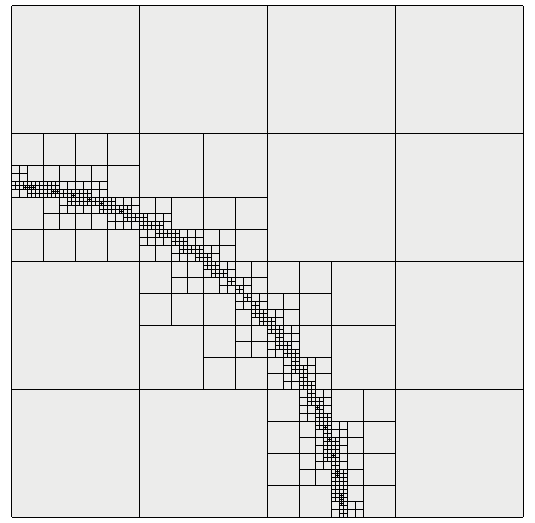} \hfill
  \includegraphics[width=0.3\textwidth]{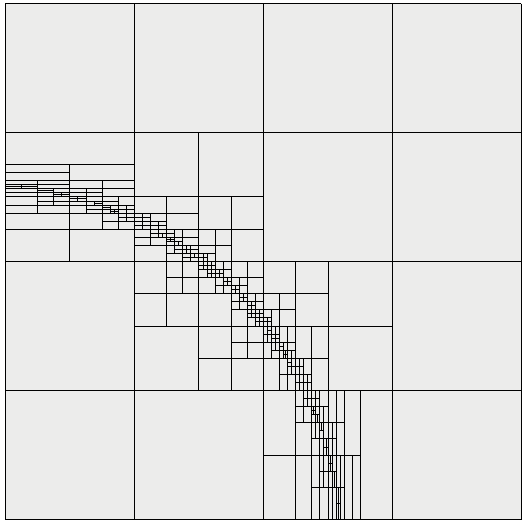}
  \caption{2D benchmark problem for finite element order $p$=2: solution (left),
    isotropic AMR mesh at 2197 DOFs (center), anisotropic AMR mesh at 1317 DOFs
    (right).}
  \label{fig:layer2d}

  \bigskip
  \includegraphics[width=\textwidth]{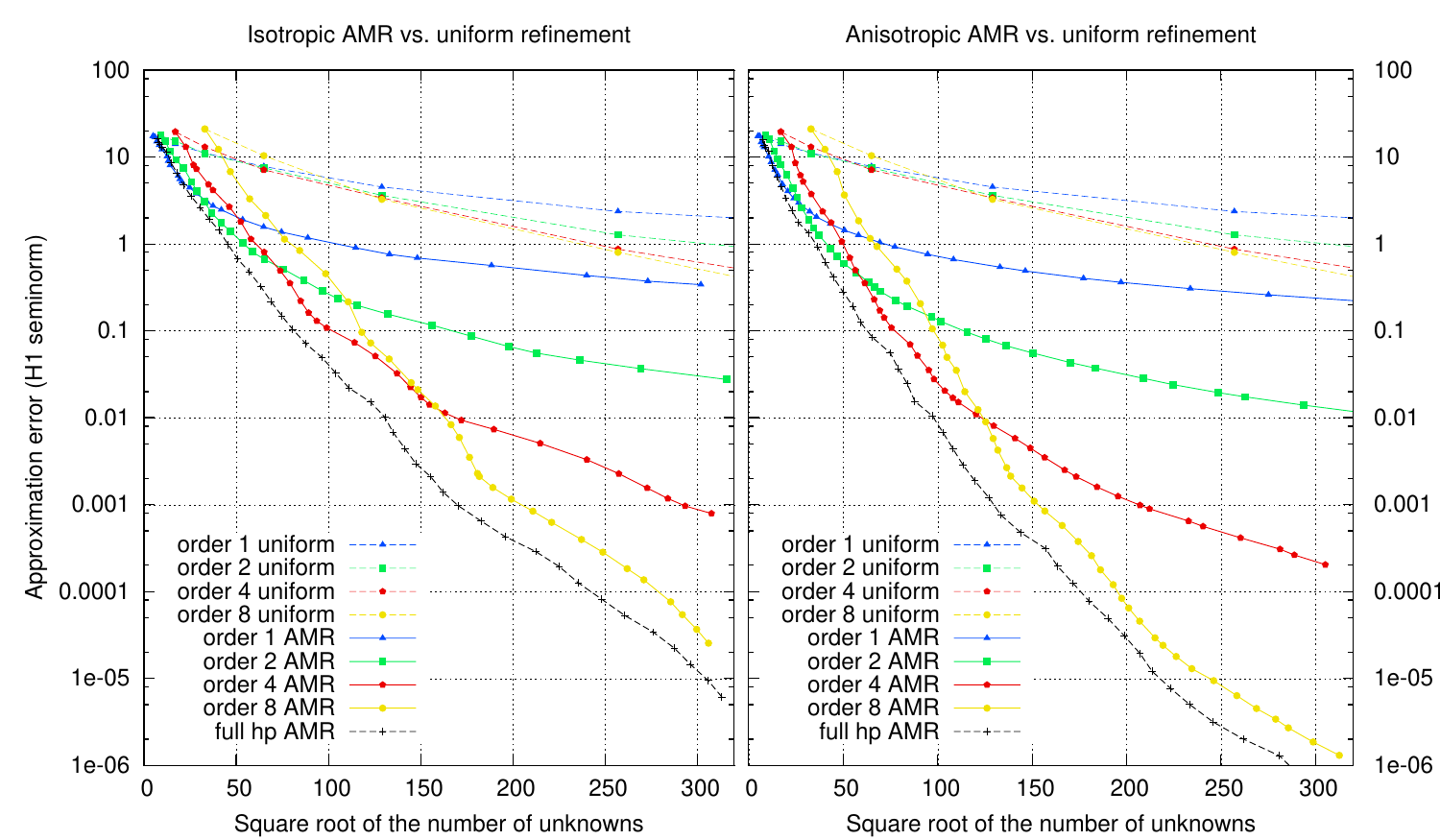}
  \caption{Convergence history for the 2D benchmark:
           isotropic refinement (left) and anisotropic refinement (right).
           The horizontal axis shows the square root of the number of
           DOFs, to make it proportional to $1/h$ of the mesh.}
  \label{fig:coverg-layer2d}
\end{figure}

\medskip
\begin{enumerate}
  \item Solve the problem on the current mesh to get an approximate solution $u_h$.
  \item For each element $K_i$, integrate the energy norm of the exact error,
        \begin{equation} e_i = ||u_h - u||_{E,K_i}. \end{equation}
  \item Record the total error $e = (\sum_i e_i^2)^{1/2}$ and the current number
        of DOFs.
  \item Refine elements for which $e_j > 0.7 \max_i\{e_i\}$.
\end{enumerate}
\medskip

Since the exact solution is not a polynomial, the errors must be calculated
with a sufficiently high order integration rule (we used a rule accurate for
order 30), especially in the first few iterations. Figure \ref{fig:layer2d}
(center) shows the mesh after 11 iterations. The convergence curves are graphed
for $p = 1, 2, 4, 8$ in Figure \ref{fig:coverg-layer2d} (left; solid lines). We
observe that the convergence histories exhibit the expected behavior.

For comparison with a previous (non-AMR) version of our code, we also plot
convergence curves for the case of uniform refinement, meaning we refine
\emph{all} elements in Step 4 above in order to avoid hanging nodes. As
expected, the benefit of local refinement is substantial. Moreover, the benefit
appears greater for higher order discretization. For completeness, we also
include the best available result for the problem (according to
\cite{Mitchell2011}), obtained with $hp$-FEM, a method
that adapts both the element size $h$ and the element order $p$ (the exact
method used is described in \cite{Cerveny12}).

\smallskip
\paragraph{Anisotropic Refinement}

Since our discretization supports anisotropic refinement, we test the same
problem again and enable splitting quadrilaterals along one axis only,
depending on the shape of the error. Assuming element $K_i$ is selected for
refinement in Step 4 based on its error $e_i$, we project and scale the
error gradient along each transformed reference space axis of the element, and
integrate the following anisotropic error indicators:
\begin{equation}
a_j = \int_{K_i} \left( J^{(j)} \cdot \nabla(u_h - u) \right)^2 \mbox{d}\mathbf{x},
\end{equation}
where $J^{(j)}$ denotes the $j$-th column of the element transformation Jacobian.
We then define the threshold $\tau_a = 0.6 / dim \, \sum a_j,$
and mark $K_i$ for refinement in its $j$-th reference axis if
$a_j > \tau_a$. Figure \ref{fig:layer2d} (right) shows the anisotropically
refined mesh after 12 iterations, and the convergence is compared to both
isotropic and uniform refinement in Figure \ref{fig:coverg-layer2d} (right).
Even though the wave front in the solution is not really aligned with the mesh,
many elements could still be refined in one direction only, which saved
up to 48\% DOFs in this problem.

\smallskip
\paragraph{3D Benchmark}

\begin{figure}
  \centering
  \includegraphics[width=0.32\textwidth]{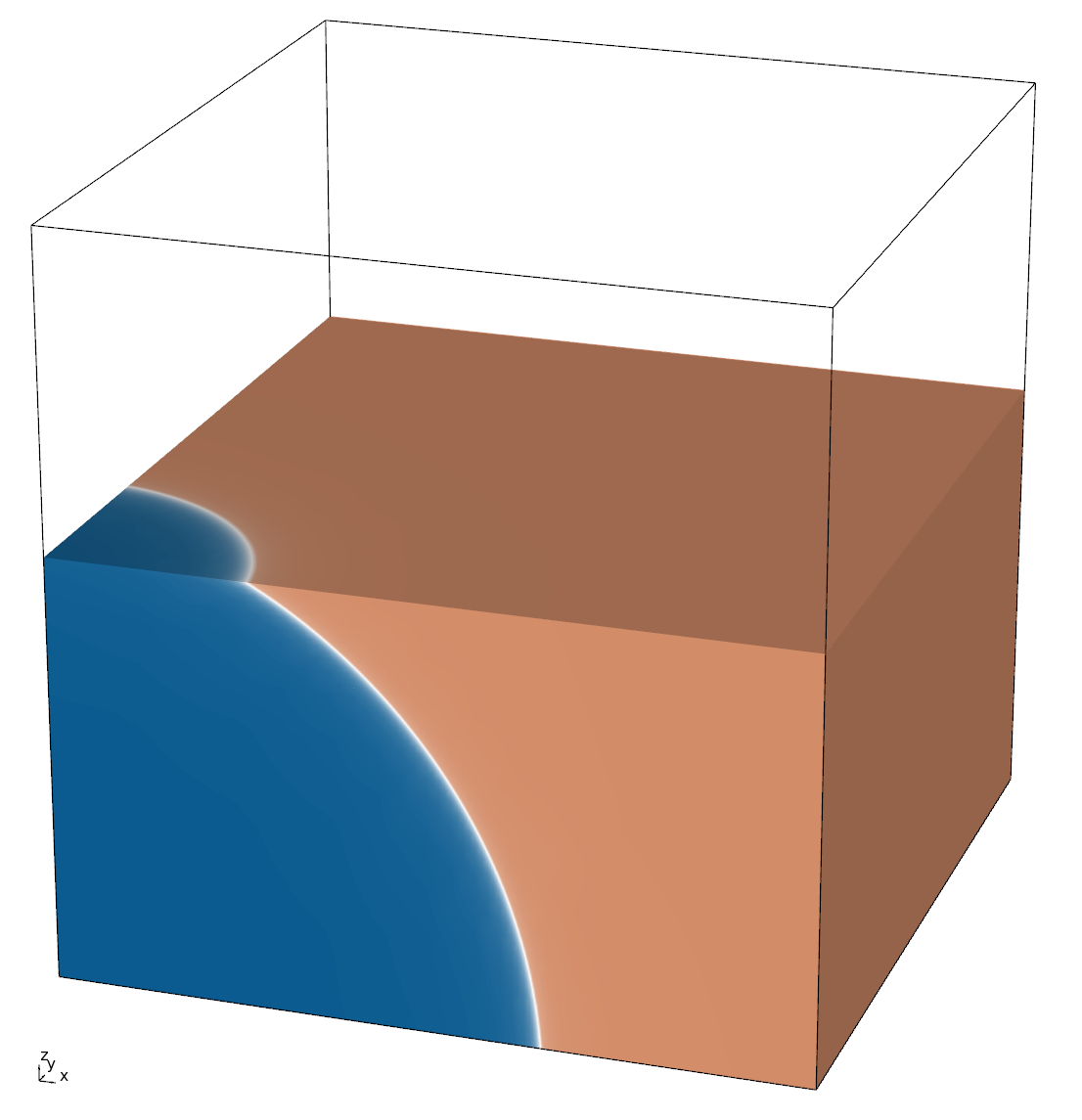} \hfill
  \includegraphics[width=0.32\textwidth]{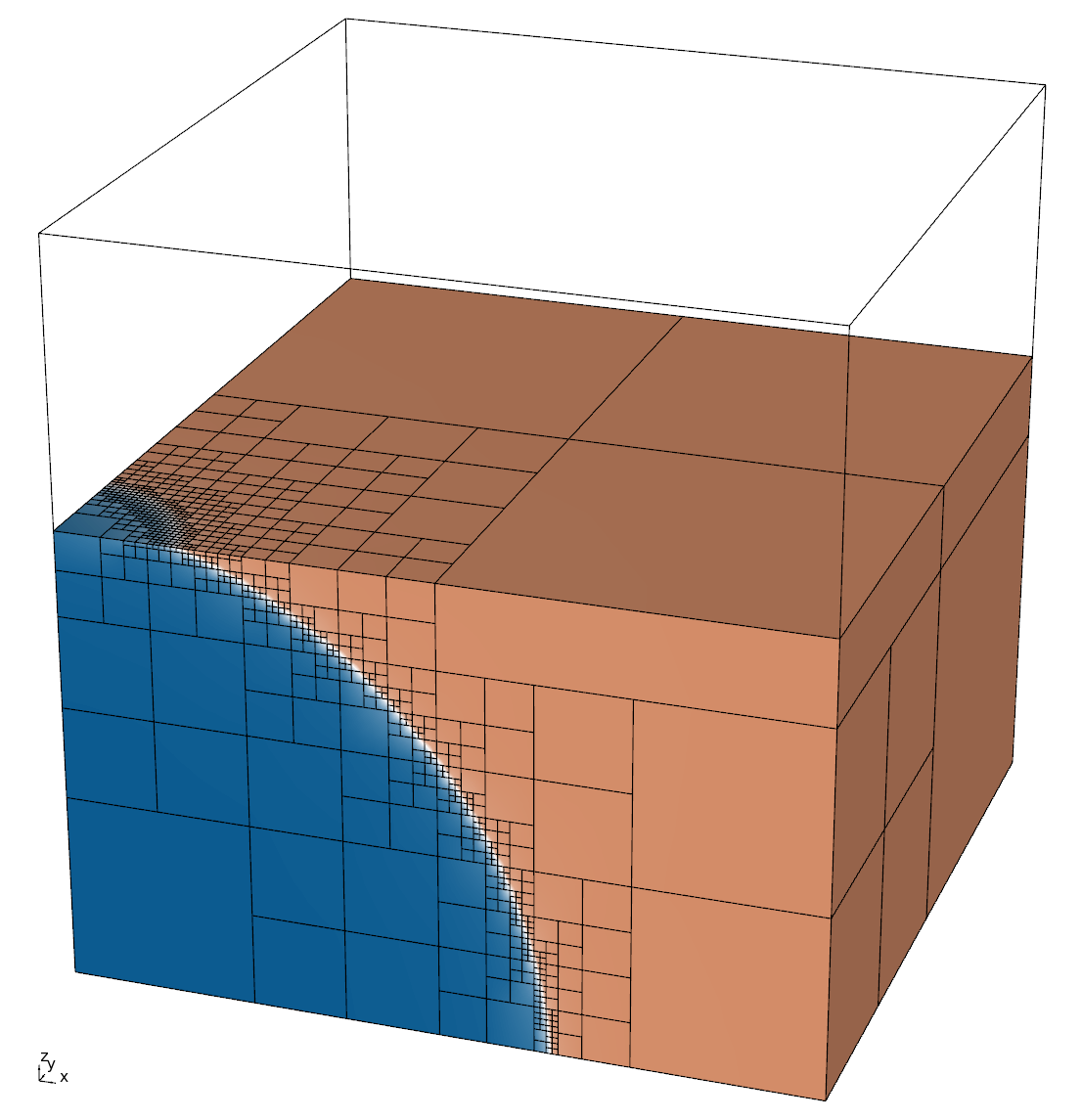} \hfill
  \includegraphics[width=0.32\textwidth]{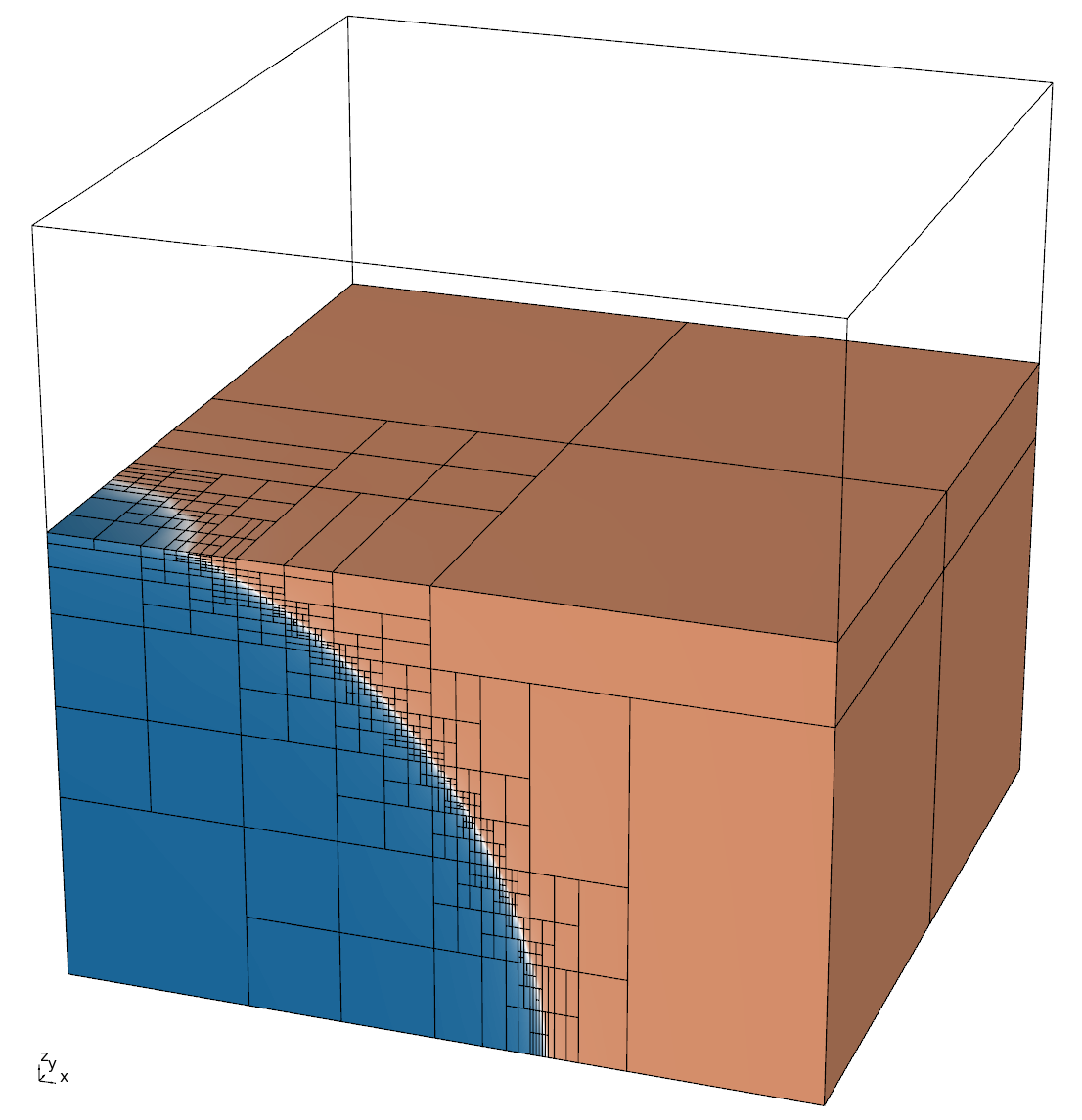}
  \caption{3D benchmark problem for order $p$=1: solution (left), isotropic AMR
    mesh at 12303 DOFs, error 5.142 (center), anisotropic AMR mesh at 5091 DOFs,
    error 4.999 (right).}
  \label{fig:layer3d}

  \bigskip
  \includegraphics[width=\textwidth]{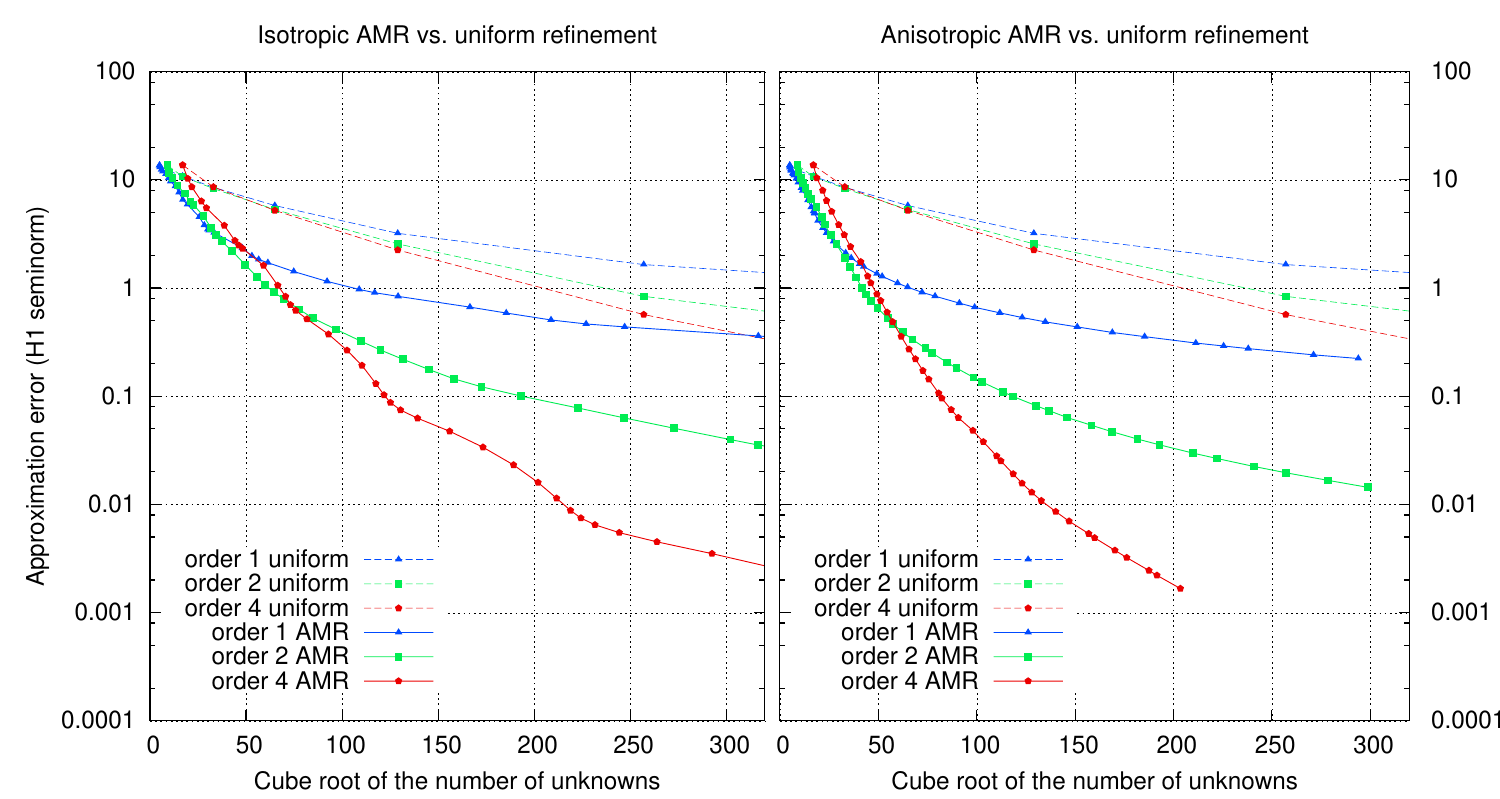}
  \caption{Convergence history for the 3D benchmark: isotropic refinement (left)
    and anisotropic refinement (right).  The horizontal axis is again normalized
    to correspond to average $1/h$ of the mesh.}
  \label{fig:converg-layer3d}
\end{figure}

Since we aim for 3D AMR, we present a straightforward generalization of the
benchmark problem into 3D. This time, the right hand side is designed for the exact
solution (shown in Figure \ref{fig:layer3d}, left) to be
\begin{equation}
u(x, y, z) = \arctan\left(\alpha\left(
   \sqrt{(x - x_c)^2 + (y - y_c)^2 + (z - z_c)^2} - r \right)\right),
\end{equation}
with similar parameters, $\alpha = 200$, $(x_c, y_c, z_c) = (-0.05, -0.05, -0.05)$
and $r = 0.7$. We start with a mesh of $4\times 4\times 4$ hexahedra
of degree $p$ and execute the same AMR loop as in 2D, for both isotropic and
anisotropic refinement.

Figure \ref{fig:layer3d} shows the isotropic and anisotropic meshes after 12 and
14 iterations, respectively. Figure \ref{fig:converg-layer3d} plots the
convergence histories of uniform refinement, isotropic AMR and anisotropic AMR.
We obtain very similar behavior as in the 2D case. The convergence curves are
smooth and do not reveal any irregularities.
The 3D version of the problem benefits even more from anisotropic refinement:
more than 80\% DOFs are saved for $p$=4. Even a slight anisotropy in the
solution means that elements with edge length ratios of 1:2, 1:4, or more, can
be used.

\smallskip
\paragraph{Curved Meshes}

One of our target applications is high-order Lagrangian hydrodynamics, where it
is routine to use high-order (curvilinear) meshes.  In MFEM, the mesh curvature
is handled simply by maintaining a~vector-valued finite element function that
represents the physical position of all nodes and, in turn, of any point within
the mesh. The associated finite element space is a subspace of $(H^1)^{dim}$,
has its own conforming prolongation matrix to keep the curvature continuous, and
does not need to coincide with the solution space.  When the mesh is refined,
the curvature function is interpolated from the original space to the finer
space. We use a version of the anisotropic 3D benchmark problem to test this
functionality on a spherical domain, as illustrated in Figure \ref{fig:curved}.

\subsection{Parallel scalability} \label{sec:scalability}

To demonstrate the parallel scalability of our algorithms, we designed a test to
put as much stress on the AMR infrastructure as possible. We choose the simple
Poisson problem again and fix the element order to $p=2$. This is the lowest
order that has both edge and face DOFs and where the stiffness matrix is still cheap
to assemble. We also omit the linear solver and instead only perform nodal
interpolation of the exact solution. This exposes the following components of
the AMR iteration:
\smallskip
\begin{itemize}
  \item constructing the $P$ matrix,
  \item assembling the parallel AMR system (using {\em hypre}'s $RAP$ function),
  \item refining elements and synchronizing ghost layers,
  \item load balancing the whole mesh at the end of the iteration.
\end{itemize}

\smallskip For the exact solution, we reuse the 3D ``wave front'' function from
the previous section. To make the problem larger, we sum two of the functions
with radii 0.2 and 0.4, and center both at $(x_c, y_c, z_c) = (0.5, 0.5, 0.5)$.
The gradient is also steeper with $\alpha = 400$. We initialize the mesh to
$32^3$ hexahedra and repeat the following steps, measuring their wall clock
times (averaged over all MPI ranks):

\smallskip
\begin{enumerate}
  \item Construct the finite element space for the current mesh (create the $P$ matrix).
  \item Assemble locally the stiffness matrix $A$ and right hand side $b$.
  \item Form the products $P^T\!A P$, $P^T b$.
  \item Eliminate Dirichlet boundary conditions from the parallel system.
  \item Project the exact solution $u$ to $u_h$ by nodal interpolation.
  \item Integrate the exact error $e_i = ||u_h - u||_{E,K_i}$ on each element.
  \item Refine elements with $e_j > 0.9 \max\{e_i\}$. Synchronize ghost layers.
  \item Load balance so each process has the same number of elements ($\pm$1):
        redistribute elements and synchronize ghost layers again.
\end{enumerate}

\begin{figure}
  \bigskip
  \centering
  \begin{minipage}{.47\textwidth}
    \centering
    \includegraphics[height=0.2\textheight]{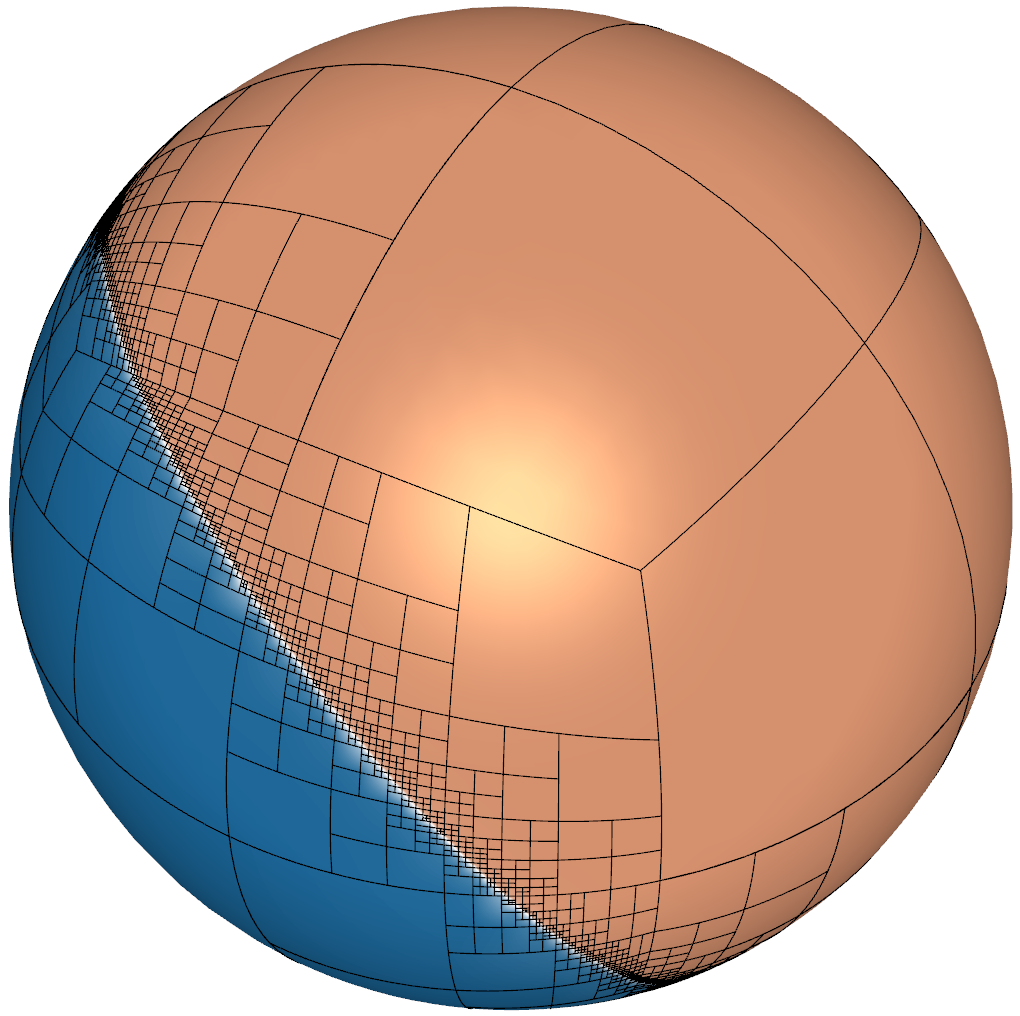}
    \caption{AMR test on a spherical domain approximated by a degree 4 curvature
             function.}
    \label{fig:curved}
  \end{minipage}
  \hfill
  \begin{minipage}{.47\textwidth}
    \centering
    \includegraphics[height=0.2\textheight]{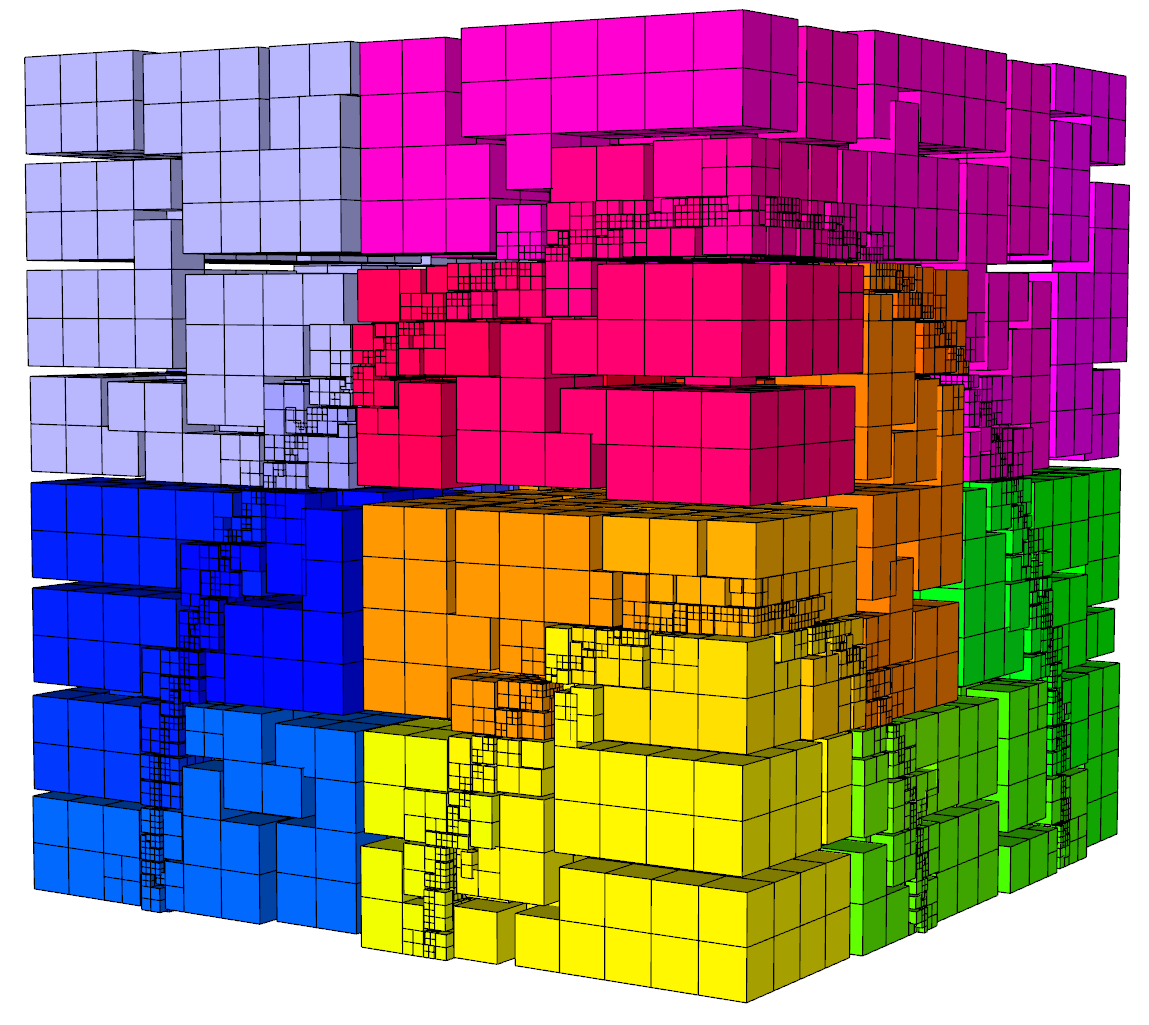}
    \includegraphics[height=0.2\textheight]{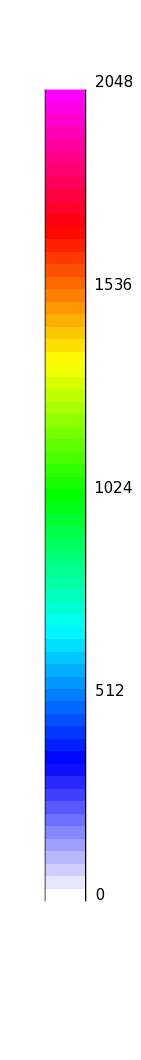}
    \caption{One octant of the parallel scaling test mesh.  Partitioning by the
      Hilbert curve is illustrated (2048 domains shown).}
    \label{fig:scaling-partition}
  \end{minipage}
\end{figure}

\smallskip Figure \ref{fig:scaling-partition} shows 1/8 of the mesh after
several iterations of this AMR loop on 16384 MPI tasks. The colors represent MPI
rank assignment.

Out of the approximately 100 iterations of the AMR loop, we select 8 iterations
that have approximately 0.5, 1, 2, 4, 8, 16, 32 and 64 million elements in the
mesh at the beginning of the iteration and plot their times as if they were 8
independent problems. This is possible because the sequence of meshes is always
the same, regardless of the number of CPU cores used. We run from 64 to 393216
(384K) cores on LLNL's Vulcan BG/Q machine. Starting with 65536 cores, there are
fewer elements in the initial mesh than parallel tasks, but this is not a
problem in our implementation.

\begin{figure}
\begin{center}
  \includegraphics[width=0.75\textwidth]{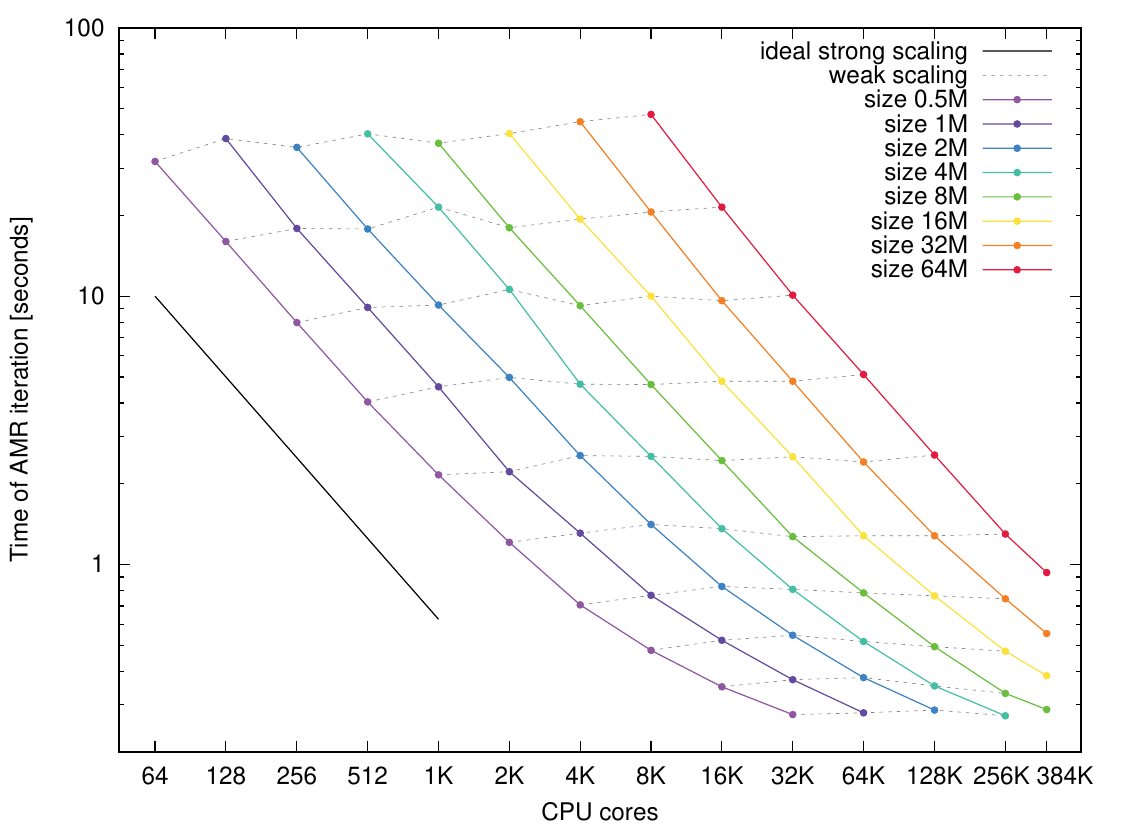}
\end{center}
\caption{Overall parallel scaling for selected iterations of the AMR loop.}
\label{fig:scaling-total}
\end{figure}

Figure \ref{fig:scaling-total} shows the total times for the selected
iterations. The solid lines show strong scaling, i.e. we follow the same AMR
iteration and its total time as the number of cores doubles. The dashed lines
skip to a twice larger problem when doubling the number of cores, showing weak
scaling, which should ideally be horizontal.

\begin{figure}
  \centering
  \includegraphics[width=0.49\textwidth]{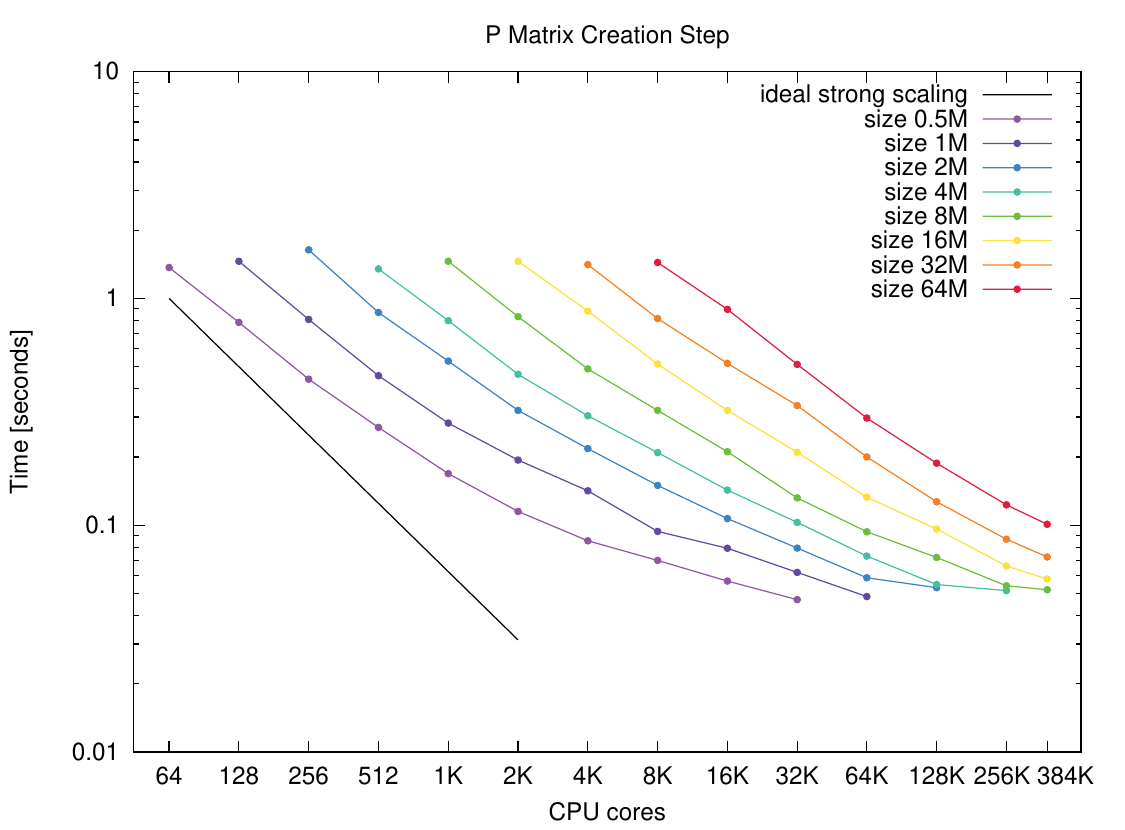}
  \includegraphics[width=0.49\textwidth]{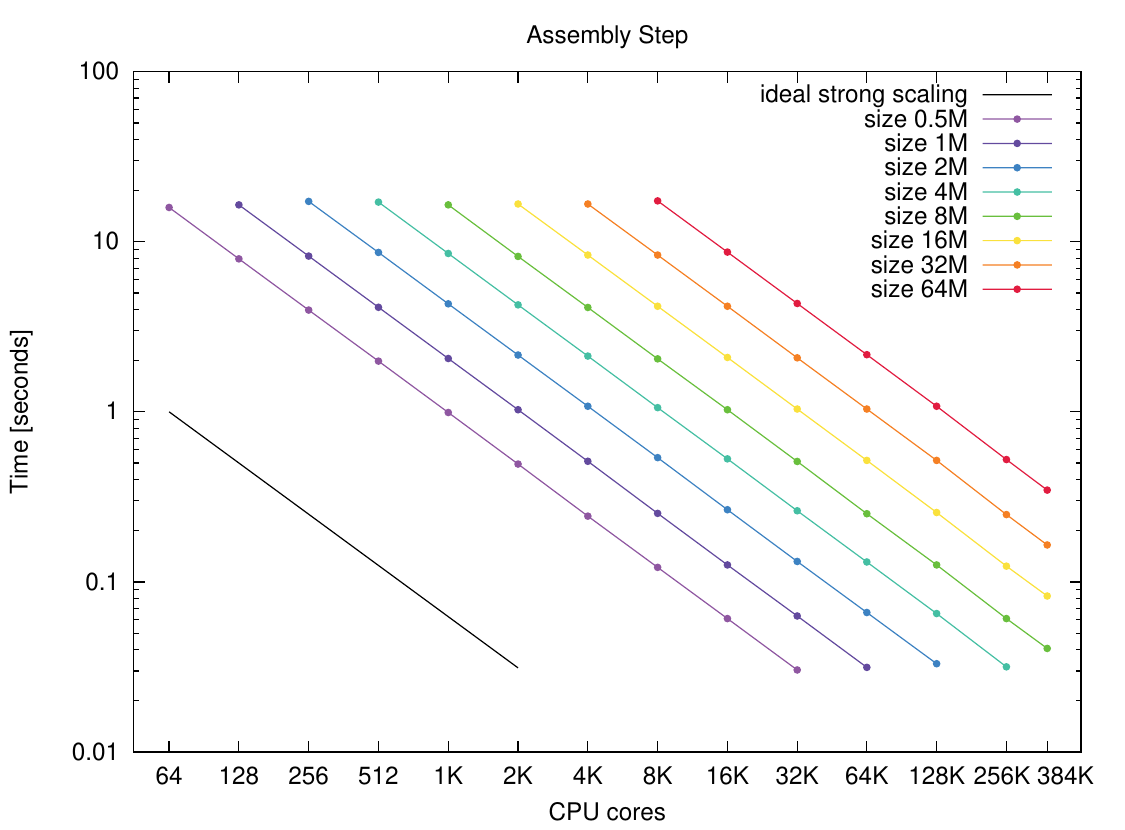}\\
  \includegraphics[width=0.49\textwidth]{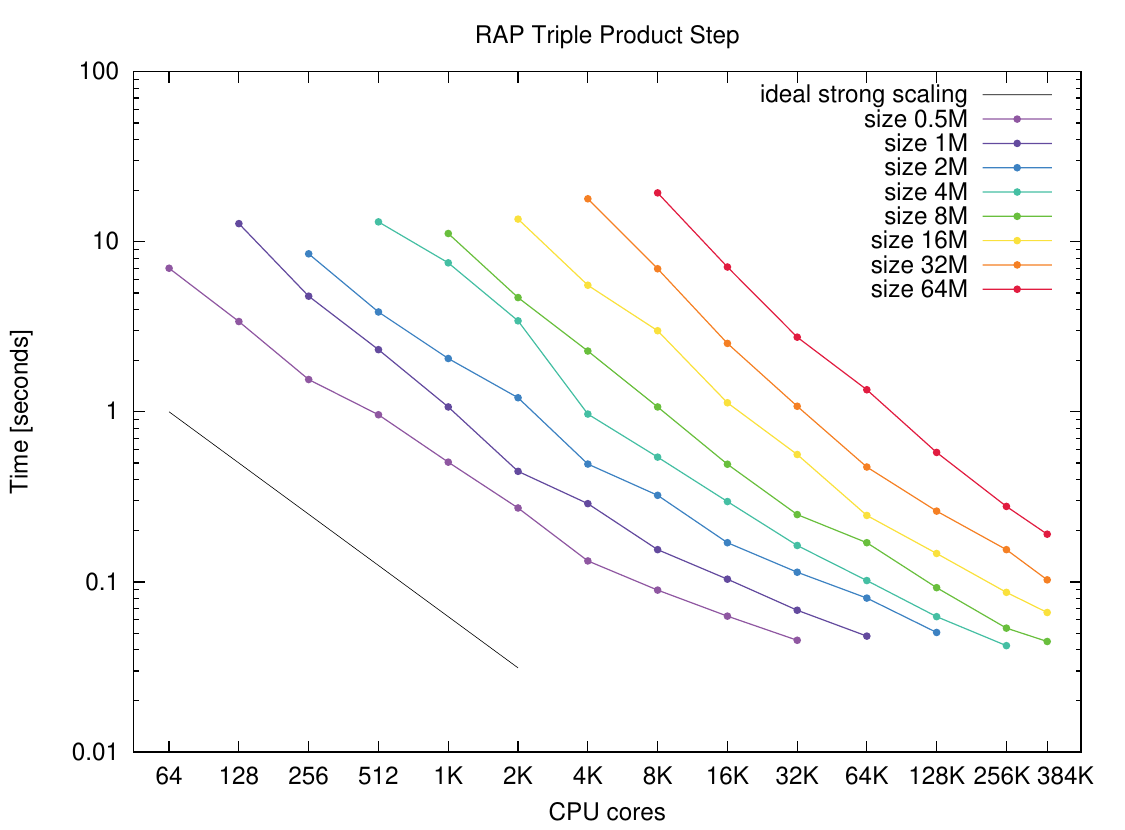}
  \includegraphics[width=0.49\textwidth]{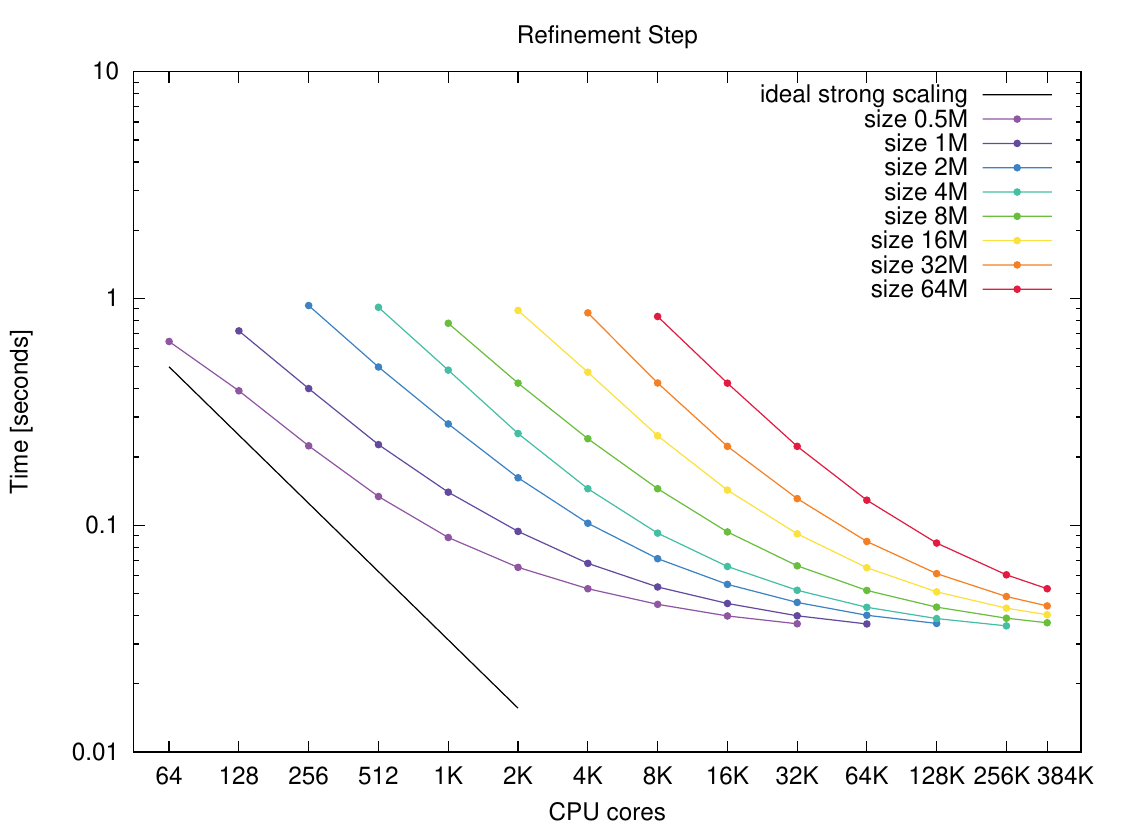}\\
  \includegraphics[width=0.49\textwidth]{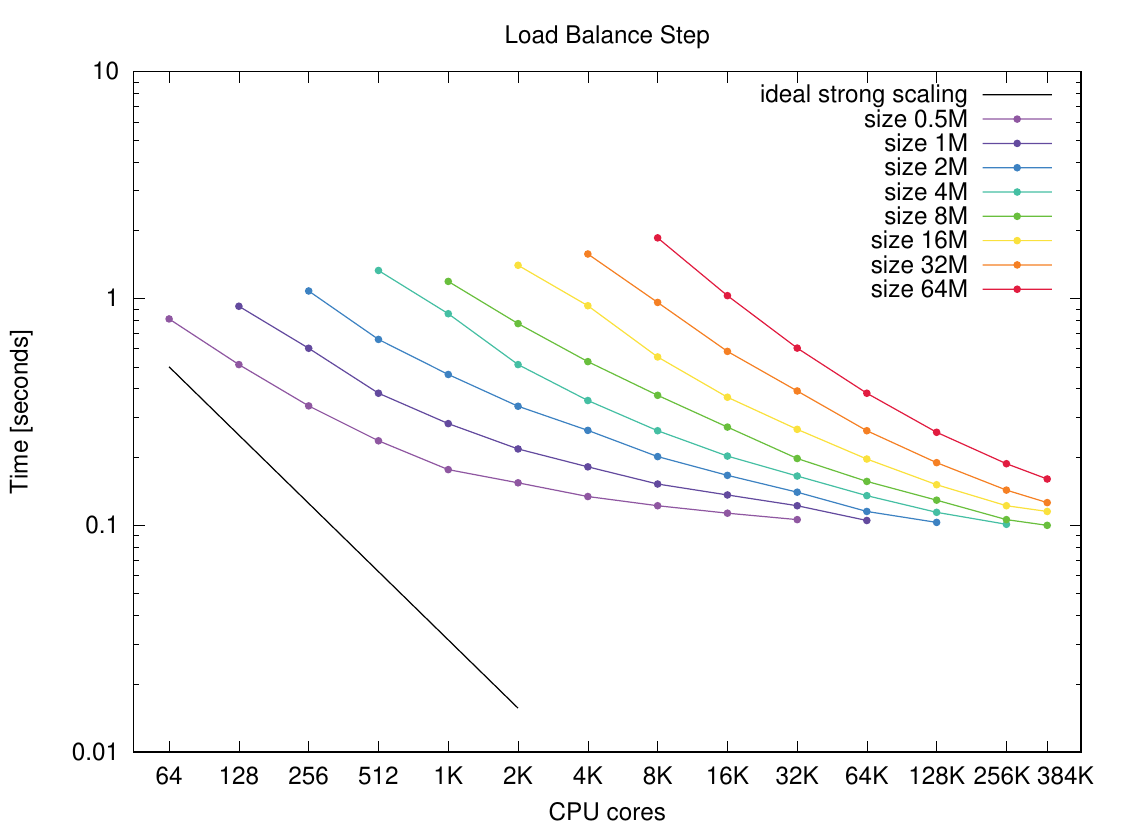}
  \includegraphics[width=0.49\textwidth]{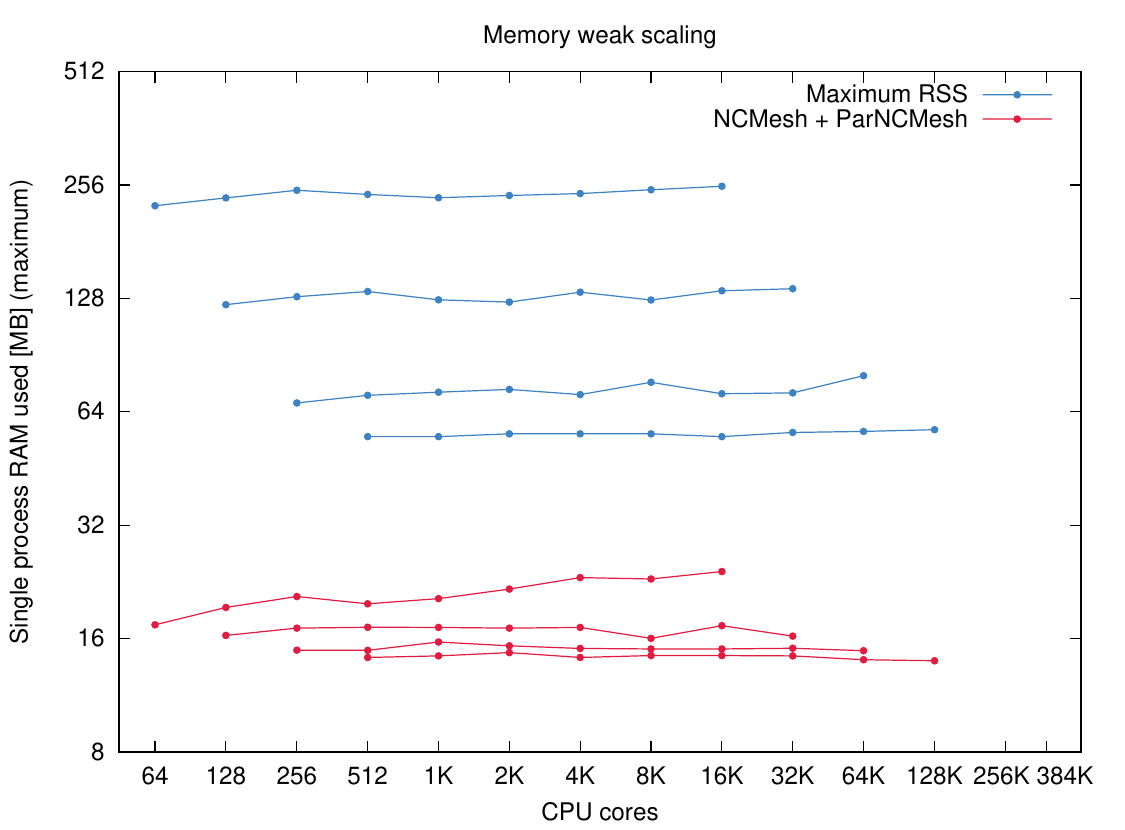}

  \caption{Scaling of individual AMR loop steps, and memory weak scaling.
           Top row: $P$ matrix construction, local assembly.
           Middle row: $RAP$ triple product step, element refinement.
           Bottom row: load balancing step, memory weak scaling.}
  \label{fig:scaling-steps}
\end{figure}

Figure \ref{fig:scaling-steps} is a partial break-down of the total iteration
time, showing individual scaling of dominant components of the iteration: $P$
matrix construction, local assembly, $RAP$ triple product, element refinement and
load balancing. Local assembly scales perfectly because it contains no
communication and the mesh is always perfectly balanced. The remaining steps do
communicate, and so scale with different degrees of success, the worst being
refinement and load balancing. Fortunately, these steps take relatively little
time compared to the total iteration time, even in this synthetic benchmark
containing no physics and solvers. Moreover, the load balancing time should be
well worth the perfect scaling of the no-communication parts.

The last graph in Figure \ref{fig:scaling-steps} shows the weak scaling of two
memory measurements. The first is the maximum RSS (resident set size) obtained
for each MPI process with \code{getrusage()}. This is the high watermark of the
RAM allocated for the process by the operating system. The second measurement is
the total memory (in MB) used by the non-conforming mesh data structure
only. This is the more difficult of the two tests, as this particular data structure
stores the ghost elements and other support data that the rest of the code does
not see and that could potentially scale badly. For both measurements, the
maximum over all MPI ranks is shown. The four curves in each set correspond to
the top four weak scaling curves in Figure \ref{fig:scaling-total}. The
remaining weak scaling curves, as well as the strong scaling curves, are not shown
here.

Finally, to test support for more than $2^{32}$ DOFs, we let the problem run on
128K cores with no limit on the number of iterations, until the available RAM
was exhausted (1GB per process in this case). Execution stopped after 131
iterations when the mesh contained $1.94\cdot 10^9$ elements, corresponding to
$14.8\cdot 10^9$ unknowns.


\subsection{Dynamic AMR in Lagrangian hydrodynamics code} \label{sec:dynamic-hydro}

In this section we show that variational restriction-based AMR can be remarkably
unintrusive when it comes to integration in a real finite element application
code. For a~demonstration, we chose the open-source high-order hydrodynamics
solver Laghos \cite{laghos}, which simulates the time-dependent Euler equations
of compressible gas dynamics in a Lagrangian frame, using high-order finite
element spatial discretization and explicit high-order time-stepping. The full
numerical method for Laghos is described in \cite{Dobrev2012}. Here we only
briefly mention that the simulation state consists of material position
$x(t)$, velocity field $v(t)$ and internal energy $e(t)$
defined over the initial domain $\Omega(t_0)$. Kinematic quantities (position,
velocity) are discretized by order $p$ continuous elements, and thermodynamic
quantities (energy) are discretized by order $(p-1)$ discontinuous elements. The
discretized material position $x(t)$ corresponds directly to the mesh
curvature function, i.e., the mesh moves and deforms together with the
material. The gas density $\rho(t)$ can be calculated at any time from the
deformation $x(t)$ and is not part of the simulation state.

To obtain the change in velocity at each time step, Laghos solves the linear
system
\begin{equation}
  M \frac{dv}{dt} = f,
\end{equation}
where $M$ is the kinematic mass matrix and $f$ contains
forces that depend on $x(t)$, $e(t)$ and the equation of state.
To handle domain decomposition in a parallel computation, the system solved
is really
\begin{equation}
  P^{T}\!M P y = P^{T} f,\ \ \
    \frac{dv}{dt} = P y,
\end{equation}
where the prolongation matrix $P$ assembles contributions from shared degrees of
freedom on parallel task boundaries (note that the product $P^{T}\!M P$
may never get evaluated explicitly in the actual implementation). Since Laghos
is already using variational restriction for domain decomposition, running the
simulation on a parallel non-conforming mesh amounts to switching the $P$ matrix
to the one constructed in Section \ref{sec:parallelP}. The only minor issue is
the order of elimination of essential boundary conditions, as discussed at the
end of Section \ref{sec:amr}.

\begin{figure}
  \centering
  \includegraphics[width=0.85\textwidth]{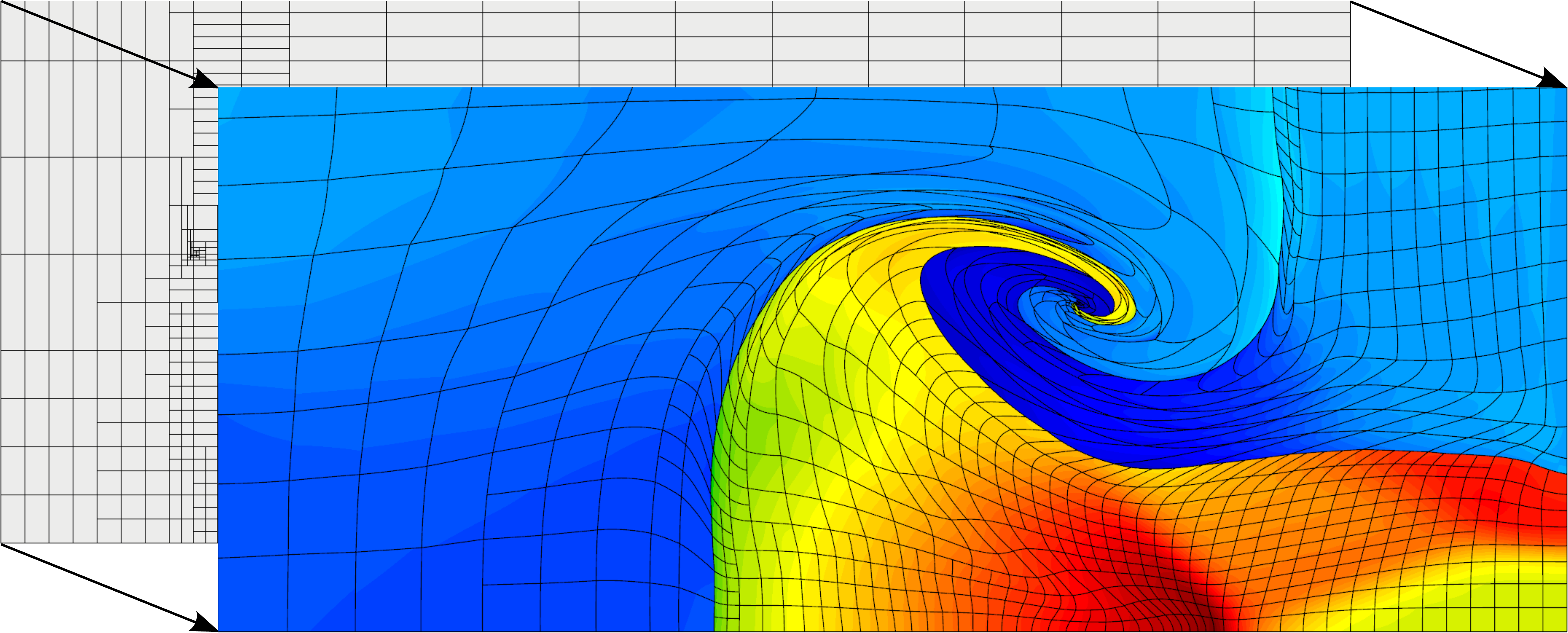}
  \caption{Static refinement in the triple-point problem.
           Initial mesh at $t = 0$ (background) refined anisotropically
           in order to obtain more regular element shapes at target time
           (foreground).}
  \label{fig:triple-pt}
\end{figure}

Our first experiment is running the 2D multi-material shock triple-point test
problem (also described in \cite{Dobrev2012}) on a static non-conforming mesh,
meaning the initial mesh is already refined and does not change during the
simulation. Based on a previous run on a coarse mesh, we anisotropically refine
elements in the three material subdomains of the initial mesh in an attempt to
counter the deformations the elements will undergo, thus improving the size of
the CFL-limited time step. The initial and target meshes are shown in Figure
\ref{fig:triple-pt}. We observe that some element deformations (diagonal
compression) cannot be countered in this way, but we conclude that the code is
capable of running on an arbitrarily refined non-conforming mesh practically
out-of-the-box.

For the second experiment, we modified the Laghos code to allow for dynamic AMR
that changes the mesh during the simulation, both via refinement and
coarsening. Changing the mesh entails interpolating the state functions on the
new mesh and reassembling both the kinematic and thermodynamic mass matrices
(which are normally constant). After each mesh change we also load balance the mesh
immediately (see Section \ref{sec:balancing}). The state vectors are migrated
together with the mesh using a parallel Boolean matrix constructed for each finite element
space as part of the load balancing operation.

Since devising a general AMR strategy for Lagrangian shock hydrodynamics is beyond the
scope of this paper, for the purpose of this demonstration we employ simple
refinement and coarsening rules tailored specifically for the well known Sedov
blast problem \cite{Dobrev2012, Sedov93}. Specifically, we look at the
artificial viscosity coefficient to detect the presence of shock discontinuity
and trigger refinement to a predefined maximum level in the undisturbed region
of the mesh ($v = 0$) in front of the shock. In the post-shock region
(nonzero~$v$) we coarsen a group of eight elements whenever their
maximum densities are below a~threshold defined as $\tau = 0.7\max\{\rho(x,
t)\},$ $x \in \Omega$. The initial mesh is the unit cube refined 5 times towards
the origin, i.e., the corner where all internal energy is deposited at time $t =
0$. This corner is never coarsened. For better stability, we also enforce
1-irregularity of the mesh at all times.

\begin{figure}
  \centering \includegraphics[width=0.315\textwidth]{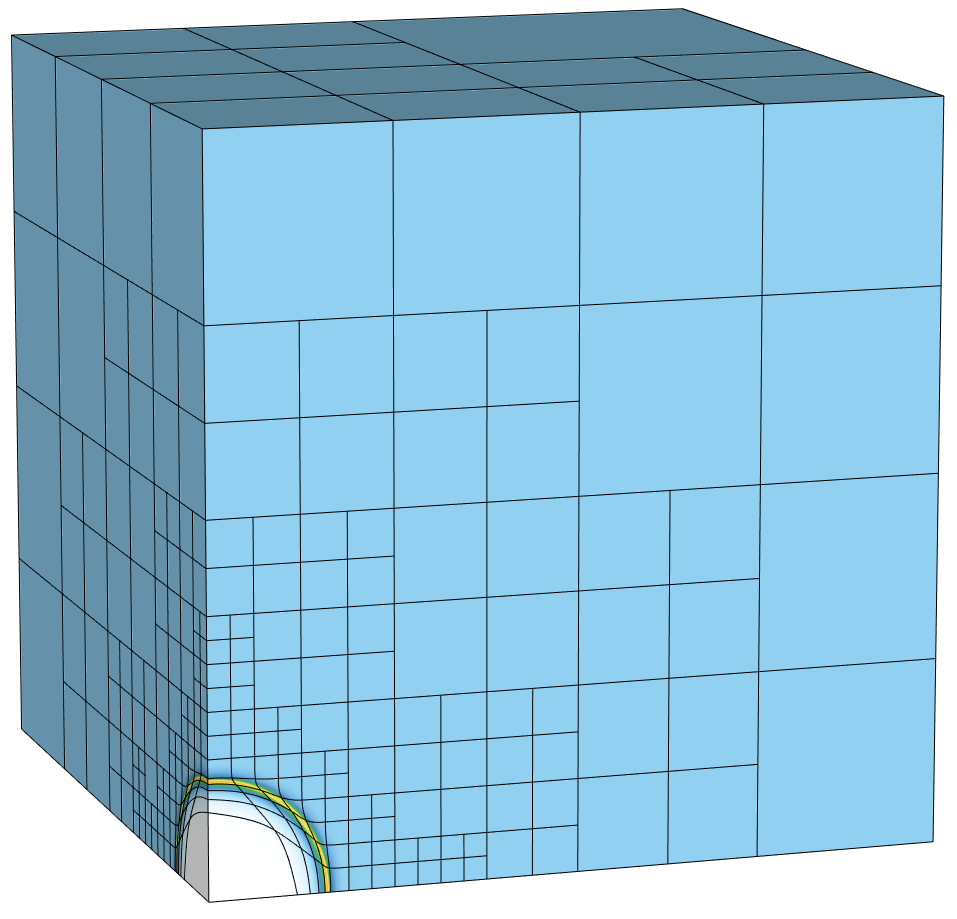}
  \hspace{1mm}
  \includegraphics[width=0.315\textwidth]{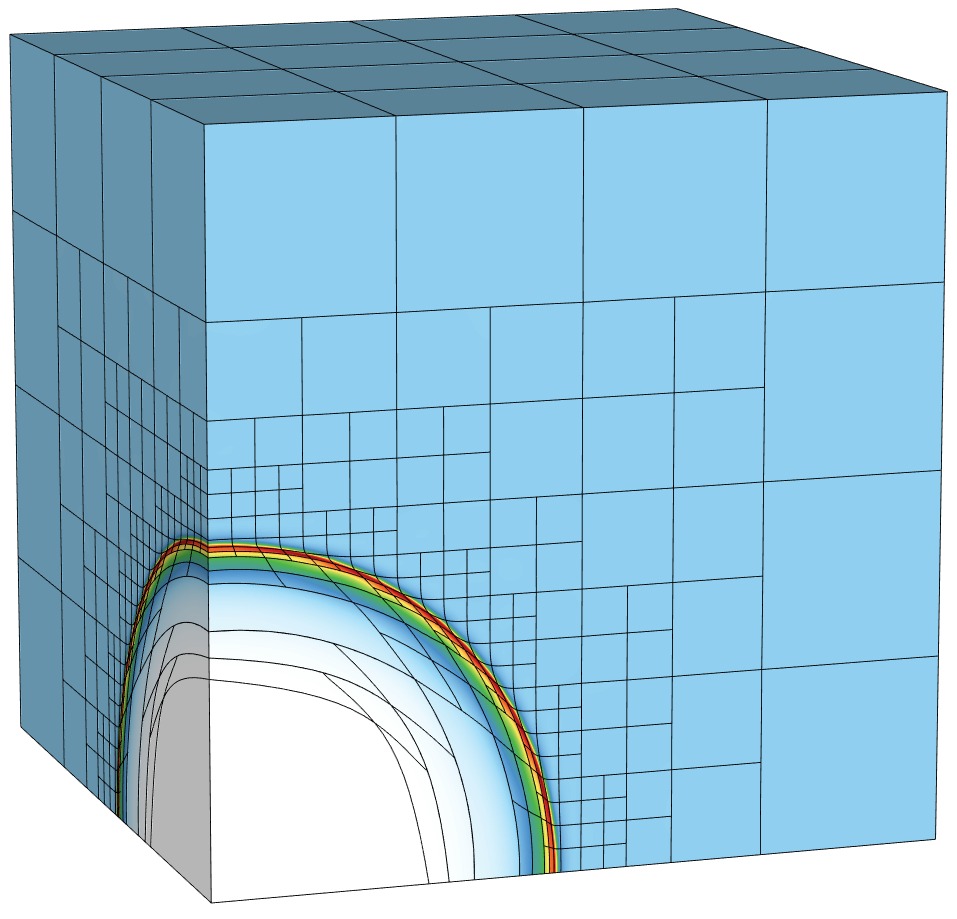}
  \hspace{1mm}
  \includegraphics[width=0.315\textwidth]{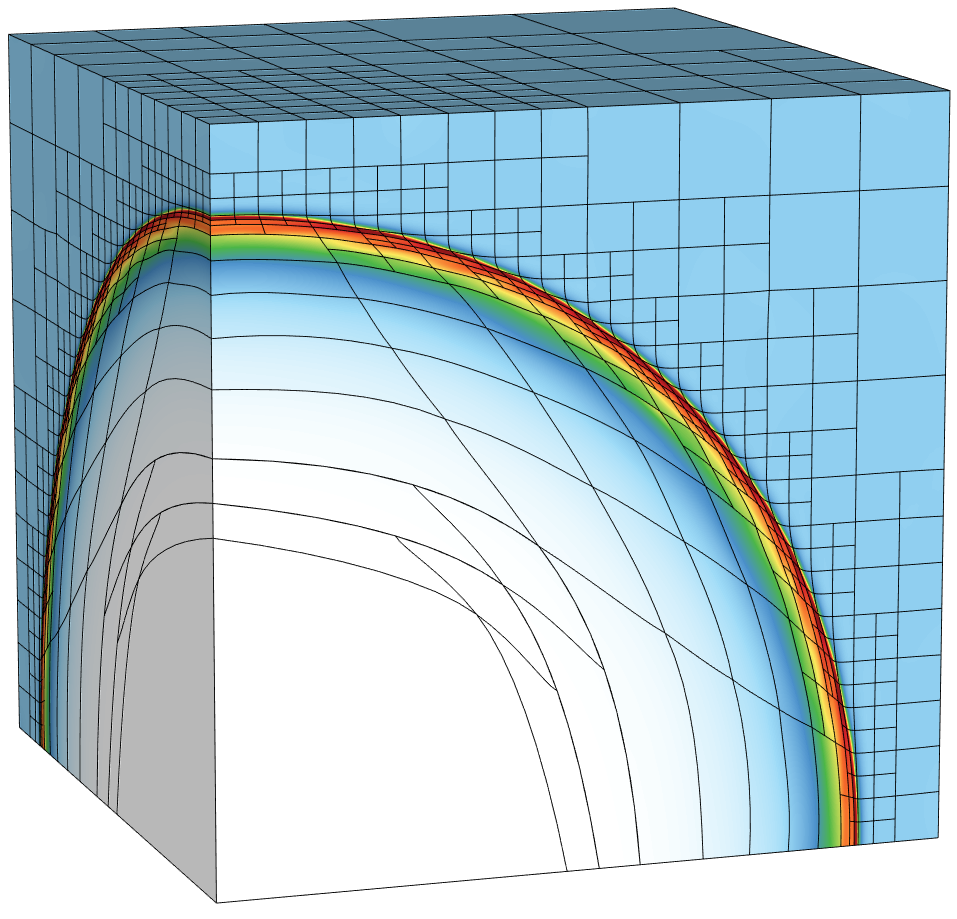}
  \caption{Dynamic refinement/coarsening in the 3D Sedov blast problem.  Mesh
    and density shown at $t = 0.0072$ (left), $t = 0.092$ (center) and $t =
    0.48$ (right). Q3Q2 elements ($p = 3$ kinematic, $p = 2$ thermodynamic
    quantities).}
  \label{fig:sedov-dynamic}
\end{figure}

With these rules we obtain a dynamic AMR computation where the maximum
refinement follows the shock, but the remaining areas are resolved with larger
elements (see Figure \ref{fig:sedov-dynamic}). This AMR-enabled version of Laghos
is publicly available in \cite{laghos} and contains
about 550 new lines of code. In contrast, obtaining these types of AMR results
in similar codes has required multiple man-years of software development in the
past.

\section{Conclusions and future work} \label{sec:conclusions}

In this paper we presented a highly scalable approach to unstructured non-conforming adaptive
mesh refinement that is easy to integrate into applications. The variational
restriction approach clearly separates the finite element assembly from the
influence of the non-conforming space. The conforming interpolation algorithm
handles a large class of triangular, quadrilateral, prismatic and hexahedral
meshes and can be used to construct arbitrary high order finite element
discretizations in $H^1$, $\Hcurl$ and $\Hdiv$, as demonstrated in examples in
this paper and in the freely available implementation in the MFEM library
\cite{mfem}.

Uniquely, our approach allows anisotropic refinement and supports arbitrarily
irregular meshes, which to the best of our knowledge has not been done
previously for continuous Galerkin hexahedral elements. Compared to purely
octree-based algorithms, we sacrifice some memory and runtime efficiency for
generality (e.g., octrees only admit isotropic refinement).  In the future, we
plan to extend the method to 4D hypercubes and perform independent space and
time adaptation, for which anisotropic refinement is crucial. We see this work
as a first step towards that goal.

\section*{Reproducibility}

Implementation of the described methods is available as open source software
\cite{mfem}. The model problems in Sections \ref{sec:model} and
\ref{sec:scalability} are similar to Example 6p in MFEM. The dynamic refinement
Sedov blast problem is available verbatim in the \code{amr} subdirectory of
the Laghos miniapp \cite{laghos}.

\bibliographystyle{siam}
\bibliography{amr}

\begin{thebibliography}{10}

\bibitem{Ainsworth97}
{\sc M.~Ainsworth and B.~Senior}, {\em Aspects of an adaptive hp-finite element
  method: Adaptive strategy, conforming approximation and efficient solvers},
  Computer Methods in Applied Mechanics and Engineering, 150 (1997),
  pp.~65--87.

\bibitem{Babuska80}
{\sc I.~Babuska and W.~Rheinboldt}, {\em Reliable error estimation and mesh
  adaptation for the finite element method},  (1979), pp.~67--108.

\bibitem{BBH2011}
{\sc W.~Bangerth, C.~Burstedde, T.~Heister, and M.~Kronbichler}, {\em
  Algorithms and data structures for massively parallel generic adaptive finite
  element codes}, ACM Trans. Math. Softw., 38 (2012), pp.~14:1--14:28.

\bibitem{dealII07}
{\sc W.~Bangerth, R.~Hartmann, and G.~Kanschat}, {\em {deal.II}---a
  general-purpose object-oriented finite element library}, ACM Trans. Math.
  Softw., 33 (2007).

\bibitem{Bank83}
{\sc R.~E. Bank, A.~H. Sherman, and A.~Weiser}, {\em Refinement algorithms and
  data structures for regular local mesh refinement}, Scientific Computing,
  (1983), pp.~3--17.

\bibitem{Zoltan12}
{\sc E.~G. Boman, U.~V. Catalyurek, C.~Chevalier, and K.~D. Devine}, {\em The
  {Z}oltan and {I}sorropia parallel toolkits for combinatorial scientific
  computing: Partitioning, ordering, and coloring}, Scientific Programming, 20
  (2012), pp.~129--150.

\bibitem{Brandt1977}
{\sc A.~Brandt}, {\em Multi-level adaptive solutions to boundary-value
  problems}, Mathematics of Computation, 31 (1977), pp.~333--390.

\bibitem{p4est11}
{\sc C.~Burstedde, L.~C. Wilcox, and O.~Ghattas}, {\em {p4est}: Scalable
  algorithms for parallel adaptive mesh refinement on forests of octrees}, SIAM
  J. Sci. Comput., 33 (2011), pp.~1103--1133.

\bibitem{Campbell2003}
{\sc P.~M. Campbell, K.~D. Devine, J.~E. Flaherty, L.~G. Gervasio, and J.~D.
  Teresco}, {\em Dynamic octree load balancing using space-filling curves},
  Tech. Rep. CS-03-01, Williams College Department of Computer Science, 2003.

\bibitem{Carey76}
{\sc G.~F. Carey}, {\em A mesh-refinement scheme for finite element
  computations}, Computer Methods in Applied Mechanics and Engineering, 7
  (1976), pp.~93 -- 105.

\bibitem{Carstensen09}
{\sc C.~Carstensen and J.~Hu}, {\em Hanging nodes in the unifying theory of a
  posteriori finite element error control}, J. Comput. Math., 27 (2009),
  pp.~215--236.

\bibitem{ceed}
{\em {CEED}: {C}enter for {E}fficient {E}xascale {D}iscretizations, {E}xascale
  {C}omputing {P}roject, {DOE}}.
\newblock \url{http://ceed.exascaleproject.org}.

\bibitem{Cerveny12}
{\sc J.~{\v{C}erven\'{y}}}, {\em Automatic hp-Adaptivity for Time-Dependent
  Problems}, PhD thesis, University of West Bohemia, 2012.

\bibitem{Ciarlet78}
{\sc P.~G. Ciarlet}, {\em The Finite Element Method for Elliptic Problems},
  North-Holland, 1978.

\bibitem{Demkowicz89}
{\sc L.~Demkowicz, J.~Oden, W.~Rachowicz, and O.~Hardy}, {\em Toward a
  universal h-p adaptive finite element strategy, part 1. constrained
  approximation and data structure}, Computer Methods in Applied Mechanics and
  Engineering, 77 (1989), pp.~79 -- 112.

\bibitem{Dobrev17}
{\sc V.~Dobrev, T.~Kolev, C.~Lee, V.~Tomov, and P.~Vassilevski}, {\em Algebraic
  hybridization and static condensation with application to scalable {H(div)}
  preconditioning}, SIAM Journal on Scientific Computing, to appear (2019).

\bibitem{Dobrev2012}
{\sc V.~A. Dobrev, T.~V. Kolev, and R.~N. Rieben}, {\em High-order curvilinear
  finite element methods for {L}agrangian hydrodynamics}, SIAM Journal on
  Scientific Computing, 34 (2012), pp.~B606--B641.

\bibitem{Flaherty97}
{\sc J.~Flaherty, R.~Loy, M.~Shephard, B.~Szymanski, J.~Teresco, and
  L.~Ziantz}, {\em Adaptive local refinement with octree load balancing for the
  parallel solution of three-dimensional conservation laws}, J. Parallel
  Distrib. Comput., 47 (1997), pp.~139--152.

\bibitem{GriebelZumbusch1997}
{\sc M.~Griebel and G.~Zumbusch}, {\em Parallel multigrid in an adaptive pde
  solver based on hashing and space-filling curves}, 1997.

\bibitem{Heuveline07}
{\sc V.~Heuveline and F.~Schieweck}, {\em H1-interpolation on quadrilateral and
  hexahedral meshes with hanging nodes}, Computing, 80 (2007), pp.~203--220.

\bibitem{hypre}
{\em hypre: A library of high performance preconditioners}.
\newblock \url{http://www.llnl.gov/CASC/hypre/}.

\bibitem{IBW2015}
{\sc T.~Isaac, C.~Burstedde, L.~Wilcox, and O.~Ghattas}, {\em Recursive
  algorithms for distributed forests of octrees}, SIAM Journal on Scientific
  Computing, 37 (2015), pp.~C497--C531.

\bibitem{DMPlexNC}
{\sc T.~Isaac and M.~G. Knepley}, {\em Support for non-conformal meshes in
  {PETSc}'s {DMPlex} interface}, CoRR, abs/1508.02470 (2015).

\bibitem{Karypis98}
{\sc G.~Karypis and V.~Kumar}, {\em A parallel algorithm for multilevel graph
  partitioning and sparse matrix ordering}, J. Parallel Distrib. Comput., 48
  (1998), pp.~71--95.

\bibitem{libMesh06}
{\sc B.~S. Kirk, J.~W. Peterson, R.~H. Stogner, and G.~F. Carey}, {\em
  {libMesh}: A {C}++ library for parallel adaptive mesh refinement/coarsening
  simulations}, Eng. with Comput., 22 (2006), pp.~237--254.

\bibitem{DMPlex}
{\sc M.~G. Knepley and D.~A. Karpeev}, {\em Mesh algorithms for {PDE with Sieve
  I:} mesh distribution}, Sci. Program., 17 (2009), pp.~215--230.

\bibitem{laghos}
{\em {Laghos}: High-order {L}agrangian hydrodynamics miniapp}.
\newblock \url{http://github.com/CEED/Laghos}.

\bibitem{Gecko}
{\sc P.~Lindstrom}, {\em The minimum edge product linear ordering problem},
  Tech. Rep. LLNL-TR-496076, Lawrence Livermore National Laboratory, Aug. 2011.

\bibitem{mfem}
{\em {MFEM}: Modular finite element methods}.
\newblock \url{http://mfem.org}.

\bibitem{Mitchell05}
{\sc W.~F. Mitchell}, {\em Hamiltonian paths through two- and three-dimensional
  grids}, J. Res. Natl. Inst. Stand. Technol., 110 (2005), pp.~127--36.

\bibitem{Mitchell07}
{\sc W.~F. Mitchell}, {\em A refinement-tree based partitioning method for
  dynamic load balancing with adaptively refined grids}, J. Parallel Distrib.
  Comput., 67 (2007), pp.~417--429.

\bibitem{Mitchell2011}
{\sc W.~F. Mitchell and M.~A. McClain}, {\em A Survey of hp-Adaptive Strategies
  for Elliptic Partial Differential Equations}, Springer Netherlands,
  Dordrecht, 2011, pp.~227--258.

\bibitem{Reps2017}
{\sc B.~Reps and T.~Weinzierl}, {\em Complex additive geometric multilevel
  solvers for helmholtz equations on spacetrees}, ACM Trans. Math. Softw., 44
  (2017), pp.~2:1--2:36.

\bibitem{Schneiders98}
{\sc R.~Schneiders}, {\em Refining quadrilateral and hexahedral element
  meshes}, Fifth International Meshing Roundtable, 1 (1998).

\bibitem{Schonfeld95}
{\sc T.~Sch{\"o}nfeld}, {\em Adaptive mesh refinement methods for
  three-dimensional inviscid flow problems}, International Journal of
  Computational Fluid Dynamics, 4 (1995), pp.~363--391.

\bibitem{Sedov93}
{\sc L.~I. Sedov}, {\em Similarity and Dimensional Methods in Mechanics}, CRC
  Press, 10th~ed., 1993.

\bibitem{Suttmeier07}
{\sc F.~T. Suttmeier}, {\em On concepts of {PDE}-software: The cellwise
  oriented approach in {DEAL}}, International Mathematical Forum, 2 (2007),
  pp.~1 -- 20.

\bibitem{TOG2005}
{\sc T.~Tu, D.~O'Hallaron, and O.~Ghattas}, {\em Scalable parallel octree
  meshing for terascale applications}, in Proc. of the 2005 ACM/IEEE Conf. on
  Supercomputing, 2005, pp.~4--.

\bibitem{SolinChapter4}
{\sc P.~\v{S}ol\'{\i}n}, {\em Partial Differential Equations and the Finite
  Element Method}, John Wiley \& Sons, 2005.

\bibitem{Solin08}
{\sc P.~\v{S}ol\'{\i}n, J.~\v{C}erven\'{y}, and I.~Dole\v{z}el}, {\em
  Arbitrary-level hanging nodes and automatic adaptivity in the hp-fem}, Math.
  Comput. Simul., 77 (2008), pp.~117--132.

\bibitem{WeinzierlMehl2011}
{\sc T.~Weinzierl and M.~Mehl}, {\em Peano--a traversal and storage scheme for
  octree-like adaptive cartesian multiscale grids}, SIAM Journal on Scientific
  Computing, 33 (2011), pp.~2732--2760.

\bibitem{wilson74}
{\sc E.~Wilson}, {\em The static condensation algorithm}, Int. J. Numer.
  Methods Eng., 8 (1974), pp.~199--203.

\end{thebibliography}

\end{document}